\documentclass[12pt,english]{article}
\usepackage[T1]{fontenc}
\usepackage[utf8]{luainputenc}
\usepackage{amssymb}

\makeatletter

\makeatletter

\makeatletter

\AtBeginDocument{\DeclareFontEncoding{T2A}{}{}}

\makeatletter

\newtheorem{theorem}{Theorem}

\newtheorem{lemma}{Lemma}

\newtheorem{proposition}[theorem]{Proposition}
\newtheorem{remark}[lemma]{Remark}

\date{}

\makeatother

\makeatother

\makeatother

\makeatother

\usepackage{babel}
\begin{document}

\title{Explicit asymptotic velocity of the boundary between particles and
antiparticles }

\author{V. A. Malyshev, A. D. Manita%
\thanks{Work of this author was supported by the Russian Foundation of Basic
Research (grants~09-01-00761 and~11-01-90421)%
}, A. A. Zamyatin}
\maketitle
\begin{abstract}
On the real line initially there are infinite number of particles
on the positive half-line., each having one of $K$ negative velocities
$v_{1}^{(+)},...,v_{K}^{(+)}$. Similarly, there are infinite number
of antiparticles on the negative half-line, each having one of $L$
positive velocities $v_{1}^{(-)},...,v_{L}^{(-)}$. Each particle
moves with constant speed, initially prescribed to it. When particle
and antiparticle collide, they both disappear. It is the only interaction
in the system. We find explicitly the large time asymptotics of $\beta(t)$
- the coordinate of the last collision before $t$ between particle
and antiparticle.

Keywords: phase boundary dynamics, random walks in cones, piece-wise
linear dynamical systems, one instrument market. 
\end{abstract}

\section{Introduction}

We consider one-dimensional dynamical model of the boundary between
two phases (particles and antiparticles, bears and bulls) where the
boundary moves due to reaction (annihilation, transaction) of pairs
of particles of different phases.

Assume that at time $t=0$ infinite number of $(+)$-particles and
$(-)$-particles are situated correspondingly on $R_{+}$ and $R_{-}$
and have one-point correlation functions 
\[
f_{+}(x,v)=\sum_{i=1}^{K}\rho_{i}^{(+)}(x)\delta(v-v_{i}^{(+)}),\qquad f_{-}(x,v)=\sum_{j=1}^{L}\rho_{j}^{(-)}(x)\delta(v-v_{j}^{(-)})
\]
 Moreover for any $i,j$ 
\[
v_{i}^{(+)}<0,\qquad v_{j}^{(-)}>0
\]
 that is two phases move towards each other. Particles of the same
phase do not see each other and move freely with the velocities prescribed
initially. The only interaction in the system is the following. When
two particles of different phases find themselves at the same point
they immediately disappear (annihilate). It follows that the phases
stay separated, and one might call any point in-between them the phase
boundary (for example it could be the point of the last collision).
Thus the boundary trajectory $\beta(t)$ is a random piece-wise constant
function of time.

The main result of the paper is the explicit formula for the asymptotic
velocity of the boundary as the function of $2(K+L)$ parameters -
densities and initial velocities. It appears to be continuous but
at some hypersurface some first derivatives in the parameters do not
exist. This kind of phase transition has very clear interpretation:
the particles with smaller activities (velocities) cease to participate
in the boundary movement - they are always behind the boundary, that
is do not influence the market price $\beta(t)$. In this paper we
consider only the case of constant densities $\rho_{i}^{(+)},\rho_{i}^{(-)}$,
that is the period of very small volatility in the market. This simplification
allows us to get explicit formulae. In \cite{MalMan} the case $K=L=1$
was considered, however with non-constant densities and random dynamics.

Main technical tool of the proof may seem surprising (and may be of
its own interest) - we reduce this infinite particle problem to the
study of a special random walk of one particle in the orthant $R_{+}^{N}$
with $N=KL$. The asymptotic behavior of this random walk is studied
using the correspondence between random walks in $R_{+}^{N}$ and
dynamical systems introduced in \cite{VM}.

The organization of the paper is the following. In section 2 we give
exact formulation of the model and of the main result. In section
3 we introduce the correspondence between infinite particle process,
random walks and dynamical systems. In sections 4 and 5 we give the
proofs.

\section{Model and the main result}

\paragraph{Initial conditions}

At time $t=0$ on the real axis there is a random configuration of
particles, consisting of $(+)$-particles and $(-)$-particles. $(+)$-particles
and $(-)$-particles differ also by the type: denote $I_{+}=\{1,2,...,K\}$
the set of types of $(+)$-particles, and $I_{-}=\{1,2,...,L\}$ -
the set of types of $(-)$-particles. Let 
\begin{equation}
0<x_{1,k}=x_{1,k}(0)<...<x_{j,k}=x_{j,k}(0)<...\label{eq:xjk}
\end{equation}
 be the initial configuration of particles of type $k\in I_{+}$,
and 
\begin{equation}
..<y_{j,i}=y_{j,i}(0)<...<y_{1,i}=y_{1,i}(0)<0\label{eq:yjk}
\end{equation}
 be the initial configuration of particles of type $i\in I_{-}$,
where the second index is the type of the particle in the configuration.
Thus all $(+)$-particles are situated on $R_{+}$ and all $(-)$-particles
on $R_{-}$. Distances between neighbor particles of the same type
are denoted by 
\begin{eqnarray}
x_{j,k}-x_{j-1,k} & = & u_{j,k}^{(+)},\quad k\in I_{+},\quad j=1,2,...\nonumber \\
y_{j-1,i}-y_{j,i} & = & u_{j,i}^{(-)},\quad i\in I_{-},\quad j=1,2,...\label{eq:distry}
\end{eqnarray}
 where we put $x_{0,k}=y_{0,i}=0$. The random configurations corresponding
to the particles of different types are assumed to be independent.
The random distances between neighbor particles of the same type are
also assumed to be independent, and moreover identically distributed,
that is random variables $u_{j,i}^{(-)},u_{j,k}^{(+)}$ are independent
and their distribution depends only on the upper and second lower
indices. Our technical assumption is that all these distributions
are absolutely continuous and have finite means. Denote $\mu_{i}^{(-)}=Eu_{j,i}^{(-)},\:\rho_{i}^{(-)}=\left(\mu_{i}^{(-)}\right)^{-1},i\in I_{-}$
, $\mu_{k}^{(+)}=Eu_{j,k}^{(+)},\:\rho_{k}^{(+)}=\left(\mu_{k}^{(+)}\right)^{-1},k\in I_{+}$.

\paragraph{Dynamics}

We assume that all $(+)$-particles of the type $k\in I_{+}$ move
in the left direction with the same constant speed $v_{k}^{(+)}$,
where $v_{1}^{(+)}<v_{2}^{(+)}<...<v_{K}^{(+)}<0$. The $(-)$-particles
of type $i\in I_{-}$ move in the right direction with the same constant
speed $v_{i}^{(-)}$, where $v_{1}^{(-)}>v_{2}^{(-)}>...>v_{L}^{(-)}>0$.
If at some time $t$ a $(+)$-particle and a $(-)$-particle are at
the same point (we call this a collision or annihilation event), then
both disappear. Collisions between particles of different phases is
the only interaction, otherwise they do not see each other. Thus,
for example, at time $t$ the $j-$th particle of type $k\in I_{+}$
could be at the point 
\[
x_{j,k}(t)=x_{j,k}(0)+v_{k}^{(+)}t
\]
 if it will not collide with some $(-)$-particle before time $t$.
Absolute continuity of the distributions of random variables $u_{j,i}^{(-)}$,$u_{j,k}^{(+)}$
guaranties that the events when more than two particles collide, have
zero probability.

We denote this infinite particle process $\mathbf{D}(t)$.

We define the boundary $\beta(t)$ between plus and minus phases to
be the coordinate of the last collision which occurred at some time
$t'<t$. For $t=0$ we put $\beta(0)=0$. Thus the trajectories of
the random process $\beta(t)$ are piecewise constant functions, we
shall assume them continuous from the left.

\paragraph{Main result}

For any pair $(J_{-},J_{+})$ of subsets $,J_{-}\subseteq I_{-},J_{+}\subseteq I_{+},$
define the numbers
\begin{equation}
V(J_{-},J_{+})=\frac{\sum_{i\in J_{-}}v_{i}^{(-)}\rho_{i}^{(-)}+\sum_{k\in J_{+}}v_{k}^{(+)}\rho_{k}^{(+)}}{\sum_{i\in J_{-}}\rho_{i}^{(-)}+\sum_{k\in J_{+}}\rho_{k}^{(+)}},V=V(I_{-},I_{+})\label{speed}
\end{equation}
 The following condition is assumed 
\begin{equation}
\{V(J_{-},J_{+})\,:\,\, J_{-}\neq\varnothing,\, J_{+}\neq\varnothing\,\}\,\cap\,\{v_{1}^{(-)},...,v_{L}^{(-)},v_{1}^{(+)},...,v_{K}^{(+)}\}=\varnothing\,.\label{condition}
\end{equation}
 If the limit $W={\displaystyle \lim_{t\to\infty}\frac{\beta(t)}{t}}$
exists a.e., we call it the asymptotic speed of the boundary. Our
main result is the explicit formula for $W$.

\begin{theorem} \label{t-osn-123-1}

The asymptotic velocity of the boundary exists and is equal to
\[
W=V(\{1,...,L_{1}\},\{1,...,K_{1}\})
\]
 where
\begin{equation}
L_{1}=\max\left\{ l\in\left\{ 1,\ldots,L\right\} :\,\, v_{l}^{(-)}>V(\{1,...,l\},I_{+})\right\} ,\label{eq:l1}
\end{equation}
 
\begin{equation}
K_{1}=\max\left\{ k\in\left\{ 1,\ldots,K\right\} :\,\, v_{k}^{(+)}<V(I_{-},\{1,...,k\})\right\} .\label{eq:k1}
\end{equation}

\end{theorem}

Note that the definition of $L_{1}$ and $K_{1}$ is not ambiguous
because $v_{1}^{(-)}>V(\{1\},I_{+})$ and $v_{1}^{(+)}<V(I_{-},\{1\})$.

Now we will explain this result in more detail. As $v_{K}^{(+)}<0<v_{L}^{(-)}$,
there can be 3 possible orderings of the numbers $v_{L}^{(-)},v_{K}^{(+)},V$: 
\begin{enumerate}
\item $v_{K}^{(+)}<V<v_{L}^{(-)}$. In this case
\[
K_{1}=K,\quad L_{1}=L,\quad W=V
\]

\item If $v_{K}^{(+)}>V$ then $V<0$ and $K_{1}<K,\quad L_{1}=L$. Moreover,
\[
W=V(\{1,...,L\},\{1,...,K_{1}\})=\min_{k\in I_{+}}V(\{1,...,L\},\{1,...,k\})\,<\, V\,<\,0
\]

\item If $v_{L}^{(-)}<V$ then $V>0$ and $K_{1}=K,\quad L_{1}<L$. Moreover,
\[
W=V(\{1,...,L_{1}\},I_{+})=\max_{l\in I_{-}}V(\{1,...,l\},I_{+})\,>\, V\,>0
\]
 The item 1 is evident. The items 2 and 3 will be explained in section
\ref{sub:Technical-lemma}. 
\end{enumerate}

\paragraph{Another scaling}

Normally the minimal difference between consecutive prices (a tick)
is very small. Moreover, one customer can have many units of the commodity.
That is why it is natural to consider the scaled densities 
\[
\rho_{j}^{(+),\epsilon}=\epsilon^{-1}\rho_{j}^{(+)},\quad\rho_{j}^{(-),\epsilon}=\epsilon^{-1}\rho_{j}^{(-)}
\]
 for some fixed constants $\rho_{j}^{(+)},\rho_{j}^{(-)}$. Then the
phase boundary trajectory $\beta^{(\epsilon)}(t)$ will depend on
$\epsilon$. The results will look even more natural. Namely, it follows
from the main theorem, that for any $t>0$ there exists the following
limit in probability 
\[
\beta(t)=\lim_{\epsilon\to0}\beta^{(\epsilon)}(t)
\]
 that is the limiting boundary trajectory.

This scaling suggests a curious interpretation of the model - the
simplest model of one instrument (for example, a stock) market. Particle
initially at $x(0)\in R_{+}$ is the seller who wants to sell his
stock for the price $x(0)$, which is higher than the existing price
$\beta(0)$. There are $K$ groups of sellers characterized by their
activity to move towards more realistic price. Similarly the $(-)$-particles
are buyers who would like to buy a stock for the price lower than
$\beta(t)$. When seller and buyer meet each other, the transaction
occurs and both leave the market. The main feature is that the traders
do not change their behavior (speeds are constant), that is in some
sense the case of zero volatility.

There are models of the market having similar type (but very different
from ours, see \cite{ContST,Parlour,rosu_1-1}). In physical literature
there are also other one-dimensional models of the boundary movement
see in \cite{Oshanin,Khorrami}.

\paragraph{Example of phase transition }

The case $K=L=1$, that is when the activities of $(+)$-particles
are the same (and similarly for $(-)$-particles), is very simple.
There is no phase transition in this case. The boundary velocity
\begin{equation}
W=\frac{v_{1}^{(+)}\rho_{1}^{(+)}+v_{1}^{(-)}\rho_{1}^{(-)}}{\rho_{1}^{(+)}+\rho_{1}^{(-)}}\label{simpleCase}
\end{equation}
 depends analytically on the activities and densities. This is very
easy to prove because the $n$-th collision time is given by the simple
formula
\begin{equation}
t_{n}=\frac{x_{n}^{(+)}(0)-x_{n}^{(-)}(0)}{-v_{1}^{(+)}+v_{1}^{(-)}}\label{eq:tn-xv}
\end{equation}
 and $n$-th collision point is given by
\begin{equation}
x_{n}^{(+)}(0)+t_{n}v_{1}^{(+)}=x_{n}^{(-)}(0)+t_{n}v_{1}^{(-)}.\label{eq:xtv}
\end{equation}

More complicated situation was considered in \cite{MalMan}. There
the movement of $(+)$-particles has random jumps in both directions
with constant drift $v_{1}^{(+)}\neq0$ (and similarly for $(-)$-particles).
In~\cite{MalMan} \emph{the order} of particles of the same type
\emph{can be changed} with time. There are no such simple formulae
as~(\ref{eq:tn-xv}) and~(\ref{eq:xtv}) in this case. The result
is however the same as in (\ref{simpleCase}).

The phase transition appears already in case when $K=2,$ $L=1$ and,
moreover, the $(-)$-particles stand still, that is $v_{1}^{(-)}=0$.
Denote $\rho_{1}^{(-)}=\rho_{0}$, $v_{i}^{(+)}=v_{i},$ $\rho_{i}^{(+)}=\rho_{i},$
$i=1,2$. Consider the function 
\[
V_{1}(v_{1},\rho_{1})=\frac{\rho_{1}v_{1}}{\rho_{0}+\rho_{1}}\,.
\]
 It is the asymptotic speed of the boundary in the system where there
is no $(+)$-particles of type 2 at all.

Then the asymptotic velocity is the function
\[
W=V(v_{1},v_{2},\rho_{1},\rho_{2})=\frac{\rho_{1}v_{1}+\rho_{2}v_{2}}{\rho_{0}+\rho_{1}+\rho_{2}}
\]
 if $v_{2}<V_{1}$ and
\[
W=V_{1}(v_{1},\rho_{1})=\frac{\rho_{1}v_{1}}{\rho_{0}+\rho_{1}}
\]
 if $v_{2}>V_{1}.$ We see that at the point $v_{2}=V_{1}$ the function
$W$ is not differentiable in $v_{2}$.

\paragraph{Balance equations - physical evidence}

Assume that the speed $w$ of the boundary is constant. Then the $(-)$-particle
will meet the boundary only if and only if $v_{i}^{(-)}>w$. Then
the mean number of $(-)$-particles of type $i$, meeting the boundary
on the time interval $(0,t)$, is $\rho_{i}^{(-)}(v_{i}^{(-)}-w)t$.
The total number of $(-)$-particles meeting the boundary during time
$t$ is
\[
\sum_{i:\, v_{i}^{(-)}>w}\rho_{i}^{(-)}(v_{i}^{(-)}-w)
\]
Similarly, the number of $(+)$-particles meeting the boundary is
\[
\sum_{j:\, v_{j}^{(+)}<w}\rho_{j}^{(+)}(w-v_{j}^{(+)})t
\]

These numbers should be equal (balance equations), and after dividing
by $t$ this gives the equation with respect to $w$
\[
\sum_{i:\, v_{i}^{(-)}>w}\rho_{i}^{(-)}(v_{i}^{(-)}-w)=\sum_{j:\, v_{j}^{(+)}<w}\rho_{j}^{(+)}(w-v_{j}^{(+)})
\]
Note that both parts are continuous in $w$. Moreover the left (right)
side is decreasing (increasing). This defines $w$ uniquely. One can
obtain the main result from this equation.

One could think that on this way one can get rigorous proof. However
it is not so easy. We develop here different techniques, that gives
much more information about the process than simple balance equations.

\section{Random walk and dynamical system in $R_{+}^{N}$}

\paragraph{Associated random walk}

One can consider the phase boundary as a special kind of server where
the customers (particles) arrive in pairs and are immediately served.
However the situation is more involved than in standard queuing theory,
because the server moves, and correlation between its movement and
arrivals is sufficiently complicated. That is why this analogy does
not help much. However we describe the crucial correspondence between
random walks in $R_{+}^{N}$ and the infinite particle problem defined
above, that allows to get the solution.

Denote $b_{i}^{(-)}(t)$ ($b_{k}^{(+)}(t)$) the coordinate of the
extreme right (left), and still existing at time~$t$, that is not
annihilated at some time $t'<t$, $(-)$-particle of type $i\in I_{-}$
($(+)$-particle of type $k\in I_{+}$). Define the distances $d_{i,k}(t)=b_{k}^{(+)}(t)-b_{i}^{(-)}(t)\geq0,\: i\in I_{-},k\in I_{+}.$
The trajectories of the random processes $b_{i}^{(-)}(t),\, b_{k}^{(+)}(t),\, d_{i,k}(t)$
are assumed left continuous. Consider the random process $D(t)=(d_{i,k}(t),(i,k)\in I)\in R_{+}^{N}$,
where $N=KL$.

Denote ${\cal \mathcal{D}}\in R_{+}^{N}$ the state space of $D(t)$.
Note that the distances $d_{i,k}(t)$, for any $t$, satisfy the following
conservation laws
\[
d_{i,k}(t)+d_{n,m}(t)=d_{i,m}(t)+d_{n,k}(t)
\]
 where $i\neq n$ and $k\neq m$. That is why the state space ${\cal \mathcal{D}}$
can be given as the set of non-negative solutions of the system of
$(L-1)(K-1)$ linear equations
\[
d_{1,1}+d_{n,m}=d_{1,m}+d_{n,1}
\]
 where $n,m\neq1$. It follows that the dimension of ${\cal \mathcal{D}}$
equals $K+L-1$. However it is convenient to speak about random walk
in $R_{+}^{N}$, taking into account that only subset of dimension
$K+L-1$ is visited by the random walk.

Now we describe the trajectories $D(t)$ in more detail. The coordinates
$d_{i,k}(t)$ decrease linearly with the speeds $v_{i}^{(-)}-v_{k}^{(+)}$
correspondingly until one of the coordinates $d_{i,k}(t)$ becomes
zero. Let $d_{i,k}(t_{0})=0$ at some time $t_{0}$. This means that
$(-)$-particle of type $i$ collided with $(+)$-particle of type
$k.$ Let them have numbers $j$ and $l$ correspondingly. Then the
components of $D(t)$ become: 
\begin{eqnarray*}
d_{i,k}(t_{0}+0) & = & u_{j+1,i}^{(-)}+u_{l+1,k}^{(+)}\\
d_{i,m}(t_{0}+0)-d_{i,m}(t_{0}) & = & u_{j+1,i}^{(-)},\,\quad m\neq k\\
d_{n,k}(t_{0}+0)-d_{n,k}(t_{0}) & = & u_{l+1,k}^{(+)},\,\quad n\neq i
\end{eqnarray*}
 and other components will not change at all, that is do not have
jumps.

Note that the increments of the coordinates $d_{n,m}(t_{0}+0)-d_{n,m}(t_{0})$
at the jump time do not depend on the history of the process before
time $t_{0}$, as the random variables. $u_{j,i}^{(-)}$($u_{j,k}^{(+)}$)
are independent and equally distributed for fixed type. It follows
that $D(t)$ is a Markov process. However that this continuous time
Markov process has singular transition probabilities (due to partly
deterministic movement). This fact however does not prevent us from
using the techniques from \cite{VM} where random walks in $Z_{+}^{N}$
were considered.

\paragraph{Ergodic case}

We call the process $D(t)$ ergodic, if there exists a neighborhood
$A$ of zero, such that the mean value $E\tau_{x}$ of the first hitting
time $\tau_{x}$ of $A$ from the point $x$ is finite for any $x\in{\cal \mathcal{D}}$.
In the ergodic case the correspondence between boundary movement and
random walks is completely described by the following theorem.

\begin{theorem}\label{thergod}

Two following two conditions are equivalent: \\
 1) The process $D(t)$ is ergodic; $\quad$ 2) $v_{K}^{(+)}<V<v_{L}^{(-)}$.

\end{theorem}

All other cases of boundary movement correspond to non-ergodic random
walks. Even more, we will see that in all other cases the process
$D(t)$ is transient. Condition~(\ref{condition}), which excludes
the set of parameters of zero measure, excludes in fact null recurrent
cases.

To understand the corresponding random walk dynamics introduce a new
family of processes.

\paragraph*{Faces}

Let $\Lambda\subseteq I=I_{-}\times I_{+}$. The face of $R_{+}^{N}$
associated with $\Lambda$ is defined as 
\begin{equation}
{\cal {\cal \mathcal{B}}}(\Lambda)=\{x\in R_{+}^{N}:\: x_{i,k}>0,\,(i,k)\in\Lambda,\, x_{i,k}=0,\,(i,k)\in\overline{\Lambda}\}\subseteq R_{+}^{N}\label{eq:B-L-def1}
\end{equation}
 If $\Lambda=\emptyset$, then ${\cal {\cal \mathcal{B}}}(\Lambda)=\{0\}$.
For shortness, instead of ${\cal {\cal \mathcal{B}}}(\Lambda)$ we
will sometimes write $\Lambda$. However, one should note that the
inclusion like $\Lambda\subset$$\Lambda_{1}$ is ALWAYS understood
for subsets of $I$, not for the faces themselves.

Define the following set of ``appropriate'' faces $\mathcal{G}=\left\{ \Lambda:\,\,\overline{\Lambda}=J_{-}\times J_{+},\,\,\, J_{-}\subseteq I_{-},\,\, J_{+}\subseteq I_{+}\right\} .$

\begin{lemma}\label{l-on-good-faces}
\[
\mathcal{D}=\bigcup_{\Lambda_{0}\in\mathcal{G}}(\mathcal{D}\cap\Lambda_{0}).
\]

\end{lemma}The proof will be given in Section 5.5. This lemma explains
why in the study of the process $D(t)$ we can consider only ``appropriate''
faces.

\subsubsection*{Induced process}

One can define a family $\mathbf{D}(t;J_{-},J_{+})$ of infinite particle
processes, where $J_{-}\subseteq I_{-},\, J_{+}\subseteq I_{+}$.
The process $\mathbf{D}(t;J_{-},J_{+})$ is the process $\mathbf{D}(t)$
with $\rho_{j}^{(+)}=0,j\notin J_{+}$ and $\rho_{j}^{(-)}=0,j\notin J_{-}$.
All other parameters (that is the densities and velocities) are the
same as for $\mathbf{D}(t)$. Note that these processes are in general
defined on different probability spaces. Obviously $\mathbf{D}(t;I_{-},I_{+})$=$\mathbf{D}(t)$.

Similarly to $\mathbf{D}(t)$, the processes $\mathbf{D}(t;J_{-},J_{+})$
have associated random walks $D(t;J_{-},J_{+})$ in $R_{+}^{N_{1}}$
with $N_{1}=|J_{-}||J_{+}|$. Usefulness of these processes is that
they describe all possible types of asymptotic behavior of the main
process $D(t)$.

Consider a face $\Lambda\in\mathcal{G}$, i.e., such face that its
complement $\overline{\Lambda}=J_{-}\times J_{+}$ where $J_{-}\subseteq I_{-}$
and $J_{+}\subseteq I_{+}$. The process $D_{\Lambda}(t)=D(t;J_{-},J_{+})=(d_{i,k}^{\Lambda}(t),(i,k)\in\overline{\Lambda})$
will be called an \textbf{\textsl{induced process}}, associated with
$\Lambda$. The coordinates $d_{i,k}^{\Lambda}(t)$ are defined in
the same way as $d_{i,k}(t)=d_{i,k}^{\Lambda}(t)$, where $\overline{\Lambda}=\{\emptyset\}$.
The state space of this process is $\mathcal{D}^{\overline{\Lambda}}=\mathcal{D}(R^{\left|\overline{\Lambda}\right|})$,
where $|\overline{\Lambda}|=|J_{-}||J_{+}|$. Face $\Lambda$ is called
\textbf{ergodic} if the induced process $D_{\Lambda}(t)$ is ergodic.

\subsubsection*{Induced vectors}

Introduce the plane
\[
\mathcal{R}(\Lambda)=\{x\in R^{N}:\: x_{i,k}=0,\,(i,k)\in\overline{\Lambda}\}\subseteq R^{N}
\]

\begin{lemma} \label{lem:iv}Let $\Lambda$ be ergodic with $\overline{\Lambda}=J_{-}\times J_{+}$,
and $D_{y}(t)$ be the process $D(t)$ with the initial point $y\in{\cal {\cal \mathcal{B}}}(\Lambda)$.
Then there exists vector $v^{\Lambda}\in\mathcal{R}(\Lambda)$ such
that for any $y\in{\cal {\cal \mathcal{B}}}(\Lambda)$ $t\geq0$,
such that $y+v^{\Lambda}t\in{\cal {\cal \mathcal{B}}}(\Lambda)$,
we have as $M\to\infty$ 
\[
\frac{D_{yM}(tM)}{M}\to y+v^{\Lambda}t
\]
 \end{lemma}

This vector $v^{\Lambda}$ will be called the \textbf{induced vector}
for the ergodic face $\Lambda$. We will see other properties of the
induced vector below.

\subsubsection*{Non-ergodic faces}

Let $\Lambda$ be the face which is not ergodic (non-ergodic face).
Ergodic face $\Lambda_{1}$: $\Lambda_{1}\supset\Lambda$ will be
called outgoing for $\Lambda$, if $v_{i,k}^{\Lambda_{1}}>0$ for
$(i,k)\in\Lambda_{1}\setminus\Lambda$. Let $\mathcal{E}(\Lambda)$
be the set of outgoing faces for the non-ergodic face $\Lambda$.

\begin{lemma}\label{nonergod}

The set $\mathcal{E}(\Lambda)$ contains the minimal element $\Lambda_{1}$
in the sense that for any $\Lambda_{2}\in\mathcal{E}(\Lambda)$ we
have $\Lambda_{2}\supseteq\Lambda_{1}.$

\end{lemma}

This lemma will be proved in section \ref{sub:nonergod}.

\subsubsection*{Dynamical system\label{dynsys}}

\label{sub:subsdyns}

We define now the piece-wise constant \textbf{vector field} $v(x)$
in $\mathcal{D}$, consisting of induced vectors, as follows: $v(x)=v^{\Lambda}$
if $x$ belongs to ergodic face $\Lambda$, and $v(x)=v^{\Lambda_{1}}$
if $x$ belongs to non-ergodic face $\Lambda$, where $\Lambda_{1}$
is the minimal element of $\mathcal{E}(\Lambda)$. Let $U^{t}$ be
the \textbf{dynamical system} corresponding to this vector field.

It follows that the trajectories $\Gamma_{x}=\Gamma_{x}(t)$ of the
dynamical system are piecewise linear. Moreover, if the trajectory
hits a non-ergodic face, it leaves it immediately. It goes with constant
speed along an ergodic face until it reaches its boundary.

We call the ergodic face $\Lambda=\mathcal{L}$ final, if either $\mathcal{L}=\emptyset$
or all coordinates of the induced vector $v^{\mathcal{L}}$ are positive.
The central statement is that the dynamical system hits the final
face, stays on it forever and goes along it to infinity, if $\mathcal{L}\neq\emptyset$.

The following theorem, together with theorem \ref{thergod}, is parallel
to theorem \ref{t-osn-123-1}. That is in all 3 cases of theorem \ref{t-osn-123-1},
theorems \ref{thergod} and \ref{thtrans} describe the properties
of the corresponding random walks in the orthant.

\begin{theorem} \null\mbox{ } \label{thtrans} 
\begin{enumerate}
\item If $D(t)$ is ergodic then the origin is the fixed point of the dynamical
system $U^{t}$. Moreover, all trajectories of the dynamical system
$U^{t}$ hit $0$. 
\item Assume $v_{K}^{(+)}>V$. Then the process $D(t)$ is transient and
there exists a unique ergodic final face $\mathcal{L}$, such that
$v_{i,k}^{\mathcal{L}}>0$ for $(i,k)\in\mathcal{L}$. This face is
\begin{eqnarray*}
\mathcal{L}(L,K_{1}) & = & \{(i,k):\,\, i=1,...,L,\,\, k=K_{1}+1,...,K\}
\end{eqnarray*}
 where $K_{1}$ is defined by (\ref{eq:k1}). Moreover, all trajectories
of the dynamical system $U^{t}$ hit $\mathcal{L}(L,K_{1})$ and stay
there forever. 
\item Assume $v_{L}^{(-)}<V$. Then the process $D(t)$ is transient and
there exists a unique ergodic final face $\mathcal{L}$, such that
$v_{i,k}^{\mathcal{L}}>0$ for $(i,k)\in\mathcal{L}$. This face is
\begin{eqnarray*}
\mathcal{L}(L_{1},K) & = & \{(i,k):\,\, i=L_{1}+1,...,L,\,\, k=1,...,K\}
\end{eqnarray*}
 where $L_{1}$ is defined by (\ref{eq:l1}). Moreover, all trajectories
of the dynamical system $U^{t}$ hit $\mathcal{L}(L_{1},K)$ and stay
there forever. 
\item For any initial point $x$ the trajectory $\Gamma_{x}(t)$ has finite
number of transitions from one face to another, until it reaches $\{0\}$
or one of the final faces. 
\end{enumerate}
\end{theorem}

This theorem will be proved in section \ref{sub:trans}.

\paragraph{Simple examples of random walks and dynamical systems}

If $K=L=1$ the process $D(t)$ is a random process on $R_{+}$. It
is deterministic on $R_{+}\setminus\{0\}$ - it moves with constant
velocity $v^{(+)}-v^{(-)}$ towards the origin. When it reaches $0$
at time $t$, it jumps backwards
\[
D(t+0)=\eta
\]
 where $\eta$ has the same distribution as $u_{1}^{(+)}+u_{1}^{(-)}$.
The dynamical system coincides with $D(t)$ inside $R_{+}$, and has
the origin as its fixed point.

If $L=1,K=2$ and moreover $v_{1}^{(-)}=0$ then the state space of
the process is $R_{+}^{2}=\{(d_{11},d_{12})\}$. Inside the quarter
plane the process is deterministic and moves with velocity $(v_{1}^{(+)},v_{2}^{(+)})$.
From any point $x$ of the boundary $d_{12}=0$ it jumps to the random
point $x+\eta_{1}$, and from any point of the boundary $d_{11}=0$
it jumps to the point $x+\eta_{2}$, where $\eta_{1},\eta_{2}$ have
the same distributions as $(u_{j,1}^{(-)},u_{j,1}^{(-)}+u_{j,2}^{(+)})$
and $(u_{j,1}^{(-)}+u_{j,1}^{(+)},u_{j,1}^{(-)})$ correspondingly.
The classification results for random walks in $Z_{+}^{2}$ can be
easily transferred to this case; the dynamical system is deterministic
and has negative components of the velocity inside $R_{+}^{2}$. When
it hits one of the axes it moves along it. The velocity is always
negative along the first axis, however along second axis it can be
either negative or positive. This is the phase transition we described
above. Correspondingly the origin is the fixed point in the first
case, and has positive value of the vector field along the second
axis, in the second case.

\section{Collisions}

\paragraph{Basic process}

Now we come back to our infinite particle process $\mathbf{D}(t)$.
The collision of particles of the types $i\in I_{-},k\in I_{+}$ we
shall call the collision of type $(i,k)$. Denote 
\[
\nu_{i,k}(T)=\#\{t:d_{i,k}(t)=0,t\in[0,T]\}
\]
 the number of collisions of type $(i,k)$ on the time interval $[0,T]$.

\begin{lemma}\label{lemma0} If the process $D(t)$ is ergodic, then
the following positive limits exist a.s. 
\begin{equation}
\pi_{i,k}=\lim_{T\to\infty}\frac{\nu_{i,k}(T)}{T}>0,\;(i,k)\in I\label{eq:ffrr}
\end{equation}
 and satisfy the following system of linear equations

\begin{equation}
v_{i}^{(-)}-v_{k}^{(+)}=\sum_{(n,m)\in I_{-}\times I_{+}}(\delta(n,i)\mu_{i}^{(-)}+\delta(m,k)\mu_{k}^{(+)})\pi_{n,m},\;(i,k)\in I\label{baleq}
\end{equation}
 \end{lemma}

Proof. Remind that the collisions can be presented as follows. If
$d_{i,k}(t_{0})=0$, then for any $n,m$
\[
d_{n,m}(t_{0}+0)-d_{n,m}(t_{0})=\delta(n,i)u_{j+1,i}^{(-)}+\delta(m,k)u_{l+1,k}^{(+)}
\]
 where $\delta(n,i)=1$ for $n=i$ and $\delta(n,i)=0$ for $n\neq i$.
Note that the proof of (\ref{eq:ffrr}) is similar to the proof of
the corresponding assertion in \cite{KMR}. For large $t$ we have
\[
d_{i,k}(t)=-(v_{i}^{(-)}-v_{k}^{(+)})t+\sum_{(n,m)\in I_{-}\times I_{+}}(\delta(n,i)\mu_{i}^{(-)}+\delta(m,k)\mu_{k}^{(+)})\nu_{n,m}(t)+o(t)
\]
 Note that this is exact equality, if instead of $\mu_{i}^{(-)}$
and $\mu_{k}^{(+)}$ we take random distances between particles. By
the law of large numbers and by (\ref{eq:ffrr}), the system (\ref{baleq})
follows.

We shall need below the following new notation. The equations (\ref{baleq})
can be rewritten in the new variables $\pi_{i}^{(-)},\pi_{k}^{(+)}$
as follows 
\[
v_{i}^{(-)}-v_{k}^{(+)}=\pi_{i}^{(-)}\mu_{i}^{(-)}+\pi_{k}^{(+)}\mu_{k}^{(+)}
\]
 where 
\[
\pi_{i}^{(-)}=\sum_{m=1}^{K}\pi_{i,m},\:\pi_{k}^{(+)}=\sum_{n=1}^{L}\pi_{n,k}
\]
 Obviously the following balance equation holds
\[
\sum_{i=1}^{L}\pi_{i}^{(-)}=\sum_{k=1}^{K}\pi_{k}^{(+)}=\sum_{i=1}^{L}\sum_{k=1}^{K}\pi_{i,k}
\]
 Rewrite the system (\ref{baleq}) in a more convenient form, using
the variables $r_{i}^{(-)}=\pi_{i}^{(-)}\mu_{i}^{(-)}$ , $r_{k}^{(+)}=\pi_{k}^{(+)}\mu_{k}^{(+)}$.
Then
\begin{eqnarray*}
v_{i}^{(-)}-v_{k}^{(+)} & = & r_{i}^{(-)}+r_{k}^{(+)},\,(i,k)\in I\\
\sum_{i=1}^{L}r_{i}^{(-)}\rho_{i}^{(-)} & = & \sum_{k=1}^{K}r_{k}^{(+)}\rho_{k}^{(+)}
\end{eqnarray*}
 It follows that for all $(i,k)\in I$
\[
v_{i}^{(-)}-r_{i}^{(-)}=r_{k}^{(+)}+v_{k}^{(+)}
\]
 Introduce the variable $w=v_{i}^{(-)}-r_{i}^{(-)}=r_{k}^{(+)}+v_{k}^{(+)}$
. We get the following system of equations with respect to the variables
$r_{i}^{(-)},\, r_{k}^{(+)},w$:

\begin{eqnarray}
v_{i}^{(-)}-r_{i}^{(-)} & = & w,\:\qquad i\in I_{-}\nonumber \\
r_{k}^{(+)}+v_{k}^{(+)} & = & w,\:\qquad k\in I_{+}\label{baleq1}\\
\sum_{i=1}^{L}r_{i}^{(-)}\rho_{i}^{(-)} & = & \sum_{k=1}^{K}r_{k}^{(+)}\rho_{k}^{(+)}\nonumber 
\end{eqnarray}
 It is easy to see that this system has the unique solution
\begin{equation}
r_{i}^{(-)}=v_{i}^{(-)}-w,\quad r_{k}^{(+)}=-v_{k}^{(+)}+w,\quad w=V\label{eq:vi}
\end{equation}
 where $V$ is defined by (\ref{speed}). If $D(t)$ is ergodic, then
by lemma \ref{lemma0} we have $r_{i}^{(-)},\, r_{k}^{(+)}>0$ for
any $i\in I_{-},\: k\in I_{+}.$

\begin{lemma}\label{lemma1}

Let the process $D(t)$ be ergodic. Then

1). $v_{K}^{(+)}<V<v_{L}^{(-)}$ .

2). The speed of the boundary $W=V$.

\end{lemma}

Proof. 1). If $D(t)$ is ergodic, then by lemma \ref{lemma0} $\pi_{i}^{(-)}>0$
and $\pi_{k}^{(+)}>0$ for all $i\in I_{-}$, $k\in I_{+}$. So, by
(\ref{eq:vi}) we have 
\[
r_{i}^{(-)}=v_{i}^{(-)}-V>0,\; r_{k}^{(+)}=-v_{k}^{(+)}+V>0
\]

2). Let $\nu_{i}^{(-)}(T)$ be the number of particles of type $i\in I_{-}$,
which had collisions during time $T$. Then 
\[
\sum_{j=1}^{\nu_{i}^{(-)}(T)}u_{j,i}^{(-)}
\]
 is the initial coordinate of the particle of type $i\in I$, which
was the last annihilated among the particle of this type. Let $T_{i}$
be the annihilation time of this particle. Then 
\[
\frac{\beta(T_{i}+0)+\sum_{j=1}^{\nu_{i}^{(-)}(T)}u_{j,i}^{(-)}}{T_{i}}=v_{i}^{(-)}
\]
 Rewrite this expression as follows
\[
\frac{\beta(T_{i}+0)-\beta(T)+\beta(T)+\sum_{j=1}^{\nu_{i}^{(-)}(T)}u_{j,i}^{(-)}}{T}=\frac{T_{i}}{T}v_{i}^{(-)}
\]
 It follows that
\[
\frac{\beta(T)}{T}=\frac{T_{i}}{T}v_{i}^{(-)}-\frac{\sum_{j=1}^{\nu_{i}^{(-)}(T)}u_{j,i}^{(-)}}{T}+\frac{\beta(T)-\beta(T_{i}+0)}{T}
\]
 By lemma \ref{lemma0} and the strong law of large numbers 
\[
\frac{\sum_{j=1}^{\nu_{i}^{(-)}(T)}u_{j,i}^{(-)}}{T}=\frac{\nu_{i}^{(-)}(T)}{T}\frac{\sum_{j=1}^{\nu_{i}^{(-)}(T)}u_{j,i}^{(-)}}{\nu_{i}^{(-)}(T)}\to\pi_{i}^{(-)}\mu_{i}^{(-)}=r_{i}^{(-)},\: a.e.
\]
 as $T\to\infty$. At the same time ergodicity of the process $D(t)$
gives that as $T\to\infty$ 
\[
\frac{T-T_{i}}{T}\to0,\;\frac{\beta(T)-\beta(T_{i}+0)}{T}\to0,\: a.e.
\]
 Thus for any $i\in I_{-}$ a.e. 
\[
\lim_{T\to\infty}\frac{\beta(T)}{T}=v_{i}^{(-)}-r_{i}^{(-)}=V
\]
 Similarly one can prove that for all $k\in I_{+}$ 
\[
\lim_{T\to\infty}\frac{\beta(T)}{T}=v_{k}^{(+)}+r_{k}^{(+)}
\]
 It follows from equations (\ref{baleq1}) and (\ref{eq:vi}) that
the boundary velocity is defined by (\ref{speed}). Lemma is proved.

\paragraph{Induced process}

Consider the faces $\Lambda$ such that $\overline{\Lambda}=J_{-}\times J_{+}$,
where $J_{-}\subseteq I_{-}$and $J_{+}\subseteq I_{+}$. Let 
\[
\nu_{i,k}^{\Lambda}(T)=\#\{t:\, d_{i,k}^{\Lambda}(t)=0,t\in[0,T]\}
\]
 be the number of collisions of type $(i,k)$ on the time interval
$[0,T]$ in the process $\mathbf{D}(t;J_{-},J_{+})$.

The following lemma is quite similar to lemma \ref{lemma0}.

\begin{lemma}If the process $D_{\Lambda}(t)$ is ergodic then the
following a.e. limits exist and are positive for all pairs $(i,k)\in\overline{\Lambda}$
\begin{equation}
\pi_{i,k}^{\Lambda}=\lim_{T\to\infty}\frac{\nu_{i,k}^{\Lambda}(T)}{T}>0\label{frec-1-1}
\end{equation}
 They satisfy the following system of linear equations 
\begin{equation}
v_{i}^{(-)}-v_{k}^{(+)}=\sum_{(n,m)\in\overline{\Lambda}}(\delta(n,i)\mu_{i}^{(-)}+\delta(m,k)\mu_{k}^{(+)})\pi_{n,m}^{\Lambda},\;(i,k)\in\overline{\Lambda}\label{baleqind}
\end{equation}
 \end{lemma}

Introduce the following notation 
\begin{eqnarray*}
\pi_{i}^{(\Lambda,-)} & = & \sum_{k\in J_{+}}\pi_{i,k}^{\Lambda},\:\quad i\in J_{-}\\
\pi_{k}^{(\Lambda,+)} & = & \sum_{i\in J_{-}}\pi_{i,k}^{\Lambda},\:\quad k\in J_{+}\\
r_{i}^{(\Lambda,-)} & = & \mu_{i}^{(-)}\pi_{i}^{(\Lambda,-)},\quad\, i\in J_{-}\\
r_{k}^{(\Lambda,+)} & = & \mu_{k}^{(+)}\pi_{k}^{(\Lambda,+)},\quad\, k\in J_{+}
\end{eqnarray*}
 For $\Lambda=\emptyset,\overline{\Lambda}=I_{-}\times I_{+}$ we
have $\pi_{i}^{(\Lambda,-)}=\pi_{i}^{(-)}$, $\pi_{k}^{(\Lambda,+)}=\pi_{k}^{(+)}$
and $r_{i}^{(\Lambda,-)}=r_{i}^{(-)}$, $r_{k}^{(\Lambda,+)}=r_{k}^{(+)}$.

Due to (\ref{baleqind}) for $(i,k)\in\overline{\Lambda}$ we have
\begin{equation}
v_{i}^{(-)}-v_{k}^{(+)}=\sum_{(n,m)\in\overline{\Lambda}}(\delta(n,i)\mu_{i}^{(-)}+\delta(m,k)\mu_{k}^{(+)})\pi_{n,m}^{\Lambda}=\mu_{i}^{(-)}\pi_{i}^{(\Lambda,-)}+\mu_{k}^{(+)}\pi_{k}^{(\Lambda,+)}=r_{i}^{(\Lambda,-)}+r_{k}^{(\Lambda,+)}\label{ind}
\end{equation}
 It follows that $v_{i}^{(-)}-r_{i}^{(\Lambda,-)}=r_{k}^{(\Lambda,+)}+v_{k}^{(+)}$
for all $(i,k)\in\overline{\Lambda}$. Put $w^{\overline{\Lambda}}=v_{i}^{(-)}-r_{i}^{(\Lambda,-)}=r_{k}^{(\Lambda,+)}+v_{k}^{(+)}$.
In this way we have obtained the following system of linear equations
(similar the system (\ref{baleq1})) with respect to variables $r_{i}^{(\Lambda,-)},r_{k}^{(\Lambda,+)},w^{\overline{\Lambda}}$:

\begin{eqnarray}
v_{i}^{(-)}-r_{i}^{(\Lambda,-)} & = & w^{\overline{\Lambda}},\quad\: i\in I_{-}\nonumber \\
r_{k}^{(\Lambda,+)}+v_{k}^{(+)} & = & w^{\overline{\Lambda}},\quad\: k\in I_{+}\label{baleq1-1}\\
\sum_{i\in J_{-}}\rho_{i}^{(-)}r_{i}^{(\Lambda,-)} & = & \sum_{k\in J_{+}}\rho_{k}^{(+)}r_{k}^{(\Lambda,+)}\nonumber 
\end{eqnarray}
 As previously, this system has the unique solution
\begin{equation}
r_{i}^{(\Lambda,-)}=v_{i}^{(-)}-w^{\overline{\Lambda}},\quad r_{k}^{(\Lambda,+)}=-v_{k}^{(+)}+w^{\overline{\Lambda}},\quad w^{\overline{\Lambda}}=V^{\overline{\Lambda}}=V(J_{-},J_{+})\label{eq:sol}
\end{equation}

For any process $\mathbf{D}(t;J_{-},J_{+})$ or for the corresponding
induced process $D_{\Lambda}(t)$(see Section 3), we also define the
boundary $\beta^{\Lambda}(t)$ as the coordinate of the last collision
$(i,k)\in\overline{\Lambda}$ before $t$. Let us assume that $\beta^{\Lambda}(0)=0$.
The trajectories of the random process $\beta^{\Lambda}(t)$ are also
piece-wise constant, we shall assume them left continuous. The following
lemma is completely analogous to lemma \ref{lemma1}.

\begin{lemma}

Let $\overline{\Lambda}=J_{-}\times J_{+}=\{i_{l},...,i_{1}\}\times\{k_{1},...,k_{m}\},$
where $i_{l}>...>i_{1}$ and $k_{1}<...<k_{m}$, and let $\Lambda$
-- be an ergodic face. Then

1). $v_{i_{l}}^{(-)}>V^{\overline{\Lambda}}=V(J_{-},J_{+})$ and $v_{k_{m}}^{(+)}<V^{\overline{\Lambda}}=V(J_{-},J_{+})$

2). The boundary velocity for the process $\mathbf{D}(t;J_{-},J_{+})$
(or for the corresponding $D_{\Lambda}(t)$) equals (with the a.e.
limit) 
\[
\lim_{t\to\infty}\frac{\beta^{\Lambda}(t)}{t}=V^{\overline{\Lambda}}=V(J_{-},J_{+})
\]

\end{lemma}

Note that $V^{\overline{\Lambda}}=V$ for $\Lambda=\emptyset$.

\begin{lemma} For any ergodic face $\Lambda$ ($\overline{\Lambda}=J_{-}\times J_{+}$)
the vector $v^{\Lambda}\in\mathcal{R}(\Lambda)$ with the coordinates
equal to
\begin{equation}
v_{i,k}^{\Lambda}=-v_{i}^{(-)}+v_{k}^{(+)}+1(i\in J_{-})\mu_{i}^{(-)}\pi_{i}^{(\Lambda,-)}+1(k\in J_{+})\mu_{k}^{(+)}\pi_{k}^{(\Lambda,+)},\quad(i,k)\in\Lambda\label{eq:iv}
\end{equation}
 is the induced vector in the sense of lemma \ref{lem:iv}.

\end{lemma}

This is quite similar to lemma 2.2, page 143 of {[}KMR{]} and lemma
4.3.2, page 87 of \cite{FMM}.

It follows from (\ref{eq:iv}) and (\ref{eq:sol}), that the coordinates
of the induced vector are given by
\begin{eqnarray}
v_{i,k}^{\Lambda} & = & -v_{i}^{(-)}+V^{\overline{\Lambda}},\,(i,k)\in\Lambda,\, i\notin J_{-},\, k\in J_{+}\label{eq:iv1}\\
v_{i,k}^{\Lambda} & =\: & v_{k}^{(+)}-V^{\overline{\Lambda}},\,(i,k)\in\Lambda,\, i\in J_{-},\, k\notin J_{+}\label{eq:iv2}\\
v_{i,k}^{\Lambda} & = & -v_{i}^{(-)}+v_{k}^{(+)},\,(i,k)\in\Lambda,\, i\notin J_{-},\, k\notin J_{+}\label{eq:iv3}\\
v_{i,k}^{\Lambda} & = & 0,\,(i,k)\in\overline{\Lambda}\nonumber 
\end{eqnarray}

Note that by condition (\ref{condition}) for all induced vectors
$v_{i,k}^{\Lambda}\neq0$ if $(i,k)\in\Lambda.$

Intuitive interpretation of this formula is the following. For example
the inequality $v_{i,k}^{\Lambda}=-v_{i}^{(-)}+V^{\overline{\Lambda}}<0,\,(i,k)\in\Lambda,\, i\notin J_{-},\, k\in J_{+}$
means that $(-)$-particles of type $i\in I_{-}$ overtake the boundary
which moves with velocity $V^{\overline{\Lambda}}$. In the contrary
case, $v_{i,k}^{\Lambda}=-v_{i}^{(-)}+V^{\overline{\Lambda}}>0$,
that is $(-)$-particles of type $i\in I_{-}$ fall behind the boundary.

\section{Proofs}

\subsection{Proof of theorem \ref{thergod} \label{sub:ergod}}

The implication $1\Rightarrow2$ has been proved in lemma \ref{lemma1}.
Now we prove that 2) implies 1). We will use the method of Lyapounov
functions to prove ergodicity. Define the Lyapounov function 
\[
f(y)=\sum_{(i,k)\in I}p_{i,k}y_{i,k}=(p,y)
\]
 where vector $p$ with coordinates $p_{i,k}>0$ will be defined below.
One has to verify the following condition: there exists $\delta>0$
such that for any ergodic face $\Lambda$, $\Lambda\neq\{0\},$
\[
f(y+v^{\Lambda})-f(y)=(p,v^{\Lambda})<-\delta
\]
 where $v^{\Lambda}$ is the induced vector corresponding to the face
$\Lambda$, see \cite{FMM}.

The system (\ref{baleq}) can be written in the matrix form 
\begin{equation}
v=A\pi\label{baleq0}
\end{equation}
 where $A$ is the $N\times N$ matrix 
\begin{equation}
A=\{a_{(i,k),(n,m)}=\delta(n,i)\mu_{i}^{(-)}+\delta(m,k)\mu_{k}^{(+)}\},\label{a}
\end{equation}
 with the elements indexed by $(i,k)\in I$, and the vector 
\begin{equation}
v=\{v_{(i,k)}=v_{i}^{(-)}-v_{k}^{(+)},(i,k)\in I\}\label{v}
\end{equation}
 It is easy to see that the coordinates of the vector $A\pi$ are
equal to
\[
(A\pi)_{i,k}=\mu_{i}^{(-)}\pi_{i}^{(-)}+\mu_{k}^{(+)}\pi_{k}^{(+)}
\]

If the assumption 2) of the theorem holds, then the system of equations
(\ref{baleq1}) has a positive solution, that is, $r_{i}^{(-)},r_{k}^{(+)}>0.$
One can choose positive $p_{i,k}$ so that the following condition
holds
\[
\pi_{i}^{(-)}=\sum_{m=1}^{K}p_{i,m},\:\pi_{k}^{(+)}=\sum_{n=1}^{L}p_{n,k}
\]
 where $\pi_{i}^{(-)}=\rho_{i}^{(-)}r_{i}^{(-)}$ and $\pi_{k}^{(+)}=\rho_{k}^{(+)}r_{k}^{(+)}$.
For example, one can put
\[
p_{i,m}=C^{-1}\pi_{i}^{(-)}\pi_{k}^{(+)}
\]
 where 
\[
C=\sum_{i=1}^{L}\pi_{i}^{(-)}=\sum_{k=1}^{K}\pi_{k}^{(+)}
\]
 Let the vector $p$ have coordinates $p_{i,k}$. Then $p$ satisfies
the system (\ref{baleq0}), that is $v=Ap$.

For ergodic face $\Lambda$ define the vector $\pi^{\Lambda}$ with
coordinates $\pi_{i,k}^{\Lambda}$, where $\pi_{i,k}^{\Lambda}$ for
$(i,k)\in\overline{\Lambda}$ are defined in (\ref{frec-1-1}) and
we put $\pi_{i,k}^{\Lambda}=0$ for $(i,k)\in\Lambda$. It follows
from (\ref{ind}) and (\ref{eq:iv}), that the induced vector can
be written as

\begin{equation}
v^{\Lambda}=-v+A\pi^{\Lambda}\label{ivektor}
\end{equation}
 with the matrix $A$ and the vector $v$ defined in (\ref{a}) and
(\ref{v}). By (\ref{ivektor}) we have
\[
v^{\Lambda}=-v+A\pi^{\Lambda}=-A(p-\pi^{\Lambda})
\]
 As the vector $A(p-\pi^{\Lambda})$ belongs to the face $\Lambda$
and $Pr_{\Lambda}\pi^{\Lambda}=0$, then 
\[
f(y+v^{\Lambda})-f(y)=(p,v^{\Lambda})=-(p,A(p-\pi^{\Lambda}))=-(p-\pi^{\Lambda},A(p-\pi^{\Lambda}))
\]

Note that the matrix $A$ in (\ref{baleq0}) is a nonnegative operator.
In fact, for any vector $y=(y_{i,j})\in R^{N}$
\begin{equation}
(Ay,y)=\sum_{i,k}(\mu_{i}^{(-)}y_{i}^{(-)}+\mu_{k}^{(+)}y_{k}^{(+)})y_{i,k}=\sum_{i=1}^{L}\mu_{i}^{(-)}\left(y_{i}^{(-)}\right)^{2}+\sum_{k=1}^{K}\mu_{k}^{(+)}\left(y_{k}^{(+)}\right)^{2}\geq0\label{pos}
\end{equation}
 where
\[
y_{i}^{(-)}=\sum_{m=1}^{K}y_{i,m},\: y_{k}^{(+)}=\sum_{n=1}^{L}y_{n,k}
\]
 Let for definiteness $\overline{\Lambda}=J_{-}\times J_{+}$. By
formula (\ref{pos})
\begin{eqnarray*}
-(p-\pi^{\Lambda},A(p-\pi^{\Lambda})) & = & -\sum_{i=1}^{L}\mu_{i}^{(-)}\left(\pi_{i}^{(-)}-\pi_{i}^{(\Lambda,-)}\right)^{2}-\sum_{k=1}^{K}\mu_{k}^{(+)}\left(\pi_{k}^{(+)}-\pi_{k}^{(\Lambda,+)}\right)^{2}\\
 & < & -\sum_{i\notin J_{-}}\mu_{i}^{(-)}\left(\pi_{i}^{(-)}\right)^{2}-\sum_{k\notin J_{+}}\mu_{k}^{(+)}\left(\pi_{k}^{(+)}\right)^{2}<0
\end{eqnarray*}
 as $\pi_{i}^{(-)},\pi_{k}^{(+)}>0$, $\pi_{i}^{(\Lambda,-)}=0$ for
$i\notin J_{-}$, $\pi_{k}^{(\Lambda,+)}=0$ if $k\notin J_{+}$.
As the number of faces is finite, one can always choose $\delta>0$,
so that 
\[
f(y+v^{\Lambda})-f(y)=-(p-\pi^{\Lambda},A(p-\pi^{\Lambda}))<-\delta
\]
 The theorem is proved.

\subsection{Proof of lemma \ref{nonergod} \label{sub:nonergod}}

For any non-ergodic face $\Lambda$ with $\overline{\Lambda}=J_{-}\times J_{+}=\{i_{1},...,i_{l}\}\times\{m_{1},...,m_{k}\},$
where $i_{1}<...<i_{l}$ and $m_{1}<...<m_{k}$, define

\[
q=\max\left\{ n\in\left\{ 1,...,l\right\} :\,\, v_{i_{n}}^{(-)}>V(\{i_{1},...,i_{n}\},\{m_{1},...,m_{k}\})\right\} ,
\]
 
\begin{equation}
r=\max\left\{ j\in\left\{ 1,...,k\right\} :\,\, v_{m_{j}}^{(+)}<V(\{i_{1},...,i_{l}\},\{m_{1},...,m_{j}\})\right\} .\label{eq:rmax}
\end{equation}
 This definition is correct because always $v_{i_{1}}^{(-)}>V(\{i_{1}\},\{m_{1},...,m_{k}\})$
and $v_{m_{1}}^{(+)}<V(\{i_{1},...,i_{l}\},\{m_{1}\})$.

Introduce the face $\Lambda_{1}$ such that $\overline{\Lambda}_{1}=\{i_{1},...,i_{q}\}\times\{m_{1},...,m_{r}\}$.
If $r=k,\; q=l$, then $v_{m_{k}}^{(+)}<V(J_{-},J_{+})<v_{i_{l}}^{(-)}$
and $\Lambda_{1}=\Lambda$. By theorem \ref{thergod} the induced
process $D_{\Lambda}(t)$ is ergodic and the face $\Lambda$ is ergodic.

So there can be two possible cases: 
\begin{itemize}
\item If $r<k,\quad q=l$, then $\overline{\Lambda}_{1}=\{i_{1},...,i_{l}\}\times\{m_{1},...,m_{r}\}$,
$v_{m_{k}}^{(+)}>V(J_{-},J_{+})$ and $V(J_{-},J_{+})<0$. 
\item If $r=k,\quad q<l$, then $\overline{\Lambda}_{1}=\{i_{1},...,i_{q}\}\times\{m_{1},...,m_{k}\}$,
$v_{i_{l}}^{(-)}<V(J_{-},J_{+})$ and$V(J_{-},J_{+})>0$. 
\end{itemize}
By construction we have $\Lambda_{1}\supset\Lambda$.

We show that $\Lambda_{1}$ is the minimal ergodic outgoing face for
$\Lambda$. Consider the first case, namely $r<k,\quad q=l$. The
second one is quite similar. Because of $v_{m_{r}}^{(+)}<V(\{i_{1},...,i_{l}\},\{m_{1},...,m_{r}\})<v_{i_{l}}^{(-)}$
we can apply theorem \ref{thergod} and so the induced process $D_{\Lambda_{1}}(t)$
is ergodic. This gives ergodicity of the face $\Lambda_{1}$.

By formula (\ref{eq:iv2}) for all $(i_{n},m_{j})\in\Lambda_{1}\setminus\Lambda=\{i_{1},...,i_{l}\}\times\{m_{r+1},...,m_{k}\}$
\[
v_{i_{n},m_{j}}^{\Lambda_{1}}=v_{m_{j}}^{(+)}-V(\{i_{1},...,i_{l}\},\{m_{1},...,m_{r}\})
\]
 and by formula (\ref{eq:rmax})
\[
v_{m_{j}}^{(+)}>V(\{i_{1},...,i_{l}\},\{m_{1},...,m_{r},m_{r+1},...,m_{j}\})
\]
 It follows from lemma \ref{tech-lem} that
\[
V(\{i_{1},...,i_{l}\},\{m_{1},...,m_{r}\})<V(\{i_{1},...,i_{l}\},\{m_{1},...,m_{r},m_{r+1},...,m_{j}\})
\]
 Thus, we get $v_{i_{n},m_{j}}^{\Lambda_{1}}>0$ for all $(i_{n},m_{j})\in\Lambda_{1}\setminus\Lambda$.
It means that the face $\Lambda_{1}$ is outgoing for $\Lambda$.

To finish the proof of lemma \ref{nonergod} it is sufficient to show
that the constructed face $\Lambda_{1}$ is the minimal outgoing face
for $\Lambda$. We give the proof by contradiction. Let there exist
an ergodic outgoing ( for $\Lambda$) face $\Lambda_{0}\supset\Lambda$
such that $\Lambda_{0}\neq\Lambda_{1}$ and $\Lambda_{1}\cap\Lambda_{0}\neq\Lambda_{1}$.
Put
\[
\overline{\Lambda}_{0}=J_{-}^{0}\times J_{+}^{0}\subset\overline{\Lambda}=\{i_{1},...,i_{l}\}\times\{m_{1},...,m_{k}\}
\]
 By (\ref{eq:iv1})-(\ref{eq:iv3}) the coordinates $v_{i,k}^{\Lambda_{0}}$
of the induced vector $v^{\Lambda_{0}}$ are given for $(i,k)\in\Lambda_{0}\setminus\Lambda$
as follows
\begin{eqnarray*}
v_{i,k}^{\Lambda_{0}} & = & -v_{i}^{(-)}+V(J_{-}^{0},J_{+}^{0}),\,(i,k)\in(J_{-}\setminus J_{-}^{0})\times J_{+}^{0},\\
v_{i,k}^{\Lambda_{0}} & = & v_{k}^{(+)}-V(J_{-}^{0},J_{+}^{0}),\,(i,k)\in J_{-}^{0}\times(J_{+}\setminus J_{+}^{0}),\,\\
v_{i,k}^{\Lambda_{0}} & = & -v_{i}^{(-)}+v_{k}^{(+)},\,(i,k)\in J_{-}\setminus J_{-}^{0}\times J_{+}\setminus J_{+}^{0}
\end{eqnarray*}
 As the face $\Lambda_{0}$ is outgoing we must have $v_{i,k}^{\Lambda_{0}}>0$
for all $(i,k)\in\Lambda_{0}\setminus\Lambda$. Thus, the only two
situations are possible: $\overline{\Lambda}_{0}=J_{-}^{0}\times\{m_{1},...,m_{k}\}$
or $\overline{\Lambda}_{0}=\{i_{1},...,i_{l}\}\times J_{+}^{0}$.
In the first case we have 
\[
v_{i,j}^{\Lambda_{0}}=-v_{i}^{(-)}+V(J_{-}^{0},\{m_{1},...,m_{k}\})>0,\,(i,j)\in(J_{-}\setminus J_{-}^{0})\times\{m_{1},...,m_{k}\}
\]
 and so $V(J_{-}^{0},\{m_{1},...,m_{k}\})>0$. But then $V(J_{-},J_{+})>0$
and this contradicts the assumption $V(J_{-},J_{+})<0$.

So $\overline{\Lambda}_{0}=\{i_{1},...,i_{l}\}\times J_{+}^{0}$.
Show that $J_{+}^{0}=\{m_{1},...,m_{r}\}$.

Let $J_{+}^{0}\neq\{m_{1},...,m_{r}\}$ and there is $j\in\{m_{1},...,m_{r}\}$
such that $j\notin J_{+}^{0}$. Then by lemma \ref{tech-lem} 
\[
v_{i,j}^{\Lambda_{0}}=v_{j}^{(+)}-V(J_{-}^{0},J_{+}^{0})<0
\]
 and, hence, the face $\Lambda_{0}$ can not be outgoing for $\Lambda$.
If $\{m_{1},...,m_{r}\}\subset J_{+}^{0}$ there exists some point
$(i_{n},m_{j})\in\Lambda_{0}\setminus\Lambda$, where $j\in\{r+1,...,k\}$
and by (\ref{eq:rmax})
\[
v_{m_{j}}^{(+)}>V(\{i_{1},...,i_{l}\},\{m_{1},...,m_{r},m_{r+1},...,m_{j}\})
\]
 It follows from theorem \ref{thergod} that the induced process $D_{\Lambda_{0}}(t)$
is non-ergodic and, hence, the face $\Lambda_{0}$ is also non-ergodic.
This contradicts the assumption on ergodicity of the face $\Lambda_{0}$.
So $J_{+}^{0}=\{m_{1},...,m_{r}\}$. Lemma is proved.

\subsection{Proof of theorem \ref{thtrans} \label{sub:trans}}

\global\long\global\long\global\long\def\gru#1#2{(#1\,|\,#2)}
 \global\long\global\long\global\long\def\tfmo{T_{f-1}}
 \global\long\global\long\def\dgm{\gamma}

The first goal of this subsection is to study trajectories $\Gamma(t)$
of the dynamical system~$U_{t}$. After that, using the obtained
knowledge about behavior of $\Gamma(t)$ we shall prove Theorem~\ref{thtrans}.
Let $\Gamma_{x}(t)$ be the trajectory of the dynamical system, starting
in the point $\Gamma_{x}(0)=x\in R_{+}^{N}$.

According to the definition of~$U_{t}$ any trajectory $\Gamma_{x}(t)$,
$t\geq0$, visits some sequence of faces. In general, this sequence
depends on the initial point~$x$ and contains ergodic and non ergodic
faces. It is very complicated to give a precise list of all faces
visited by the concrete trajectory started from a given point~$x$.
Our idea is to find a common \emph{finite subsequence} $\Lambda_{1},\Lambda_{2},...,\Lambda_{n}$
of ergodic faces in the order they are visited by any trajectory.
We find this subsequence together with the time moments $t_{1}$,
$t_{2},...,t_{n}$, where $t_{k}$ is the first time the trajectory
enters the closure of $\Lambda_{k}$. Moreover, it will follow from
our proof that the intervals $t_{k}-t_{k-1}$ are finite, the dimensions
of the ergodic faces in this sequence decrease and any trajectory,
after hitting the closure of some face in this sequence, will never
leave this closure.

\noindent \begin{proposition} \label{p-put-dyn-s} There exists a
monotone sequence of faces 
\[
\Lambda_{1}\supset\Lambda_{2}\supset\cdots\supset\Lambda_{r}\supset\cdots\supset\Lambda_{n},\qquad\dim F_{i}>\dim F_{i+1},
\]
 and a sequence of time moments 
\[
t_{1}\leq t_{2}\leq\cdots\leq t_{r}\leq\cdots\leq t_{n}<+\infty,
\]
 depending on $x$, and having the following property
\[
\Gamma_{x}(t)\in F_{r}\qquad\forall t\geq t_{r}\,,
\]
 where $F_{r}=cl(\Lambda_{r})$ denotes the closure of $\Lambda_{r}$
in $R_{+}^{N}$. Moreover, the sequence $\Lambda_{1},\Lambda_{2},...,\Lambda_{n}$
depends only on the parameters of the model (that is on the velocities
and densities), but the sequence of time moments $t_{1}$, $t_{2},...,t_{n}$
depends also on the initial point~$x$ of the trajectory~$\Gamma_{x}(t)$.
Thus any trajectory will hit the final set~$F_{fin}=F_{n}$ in finite
time.

\end{proposition}

The proof of Proposition~\ref{p-put-dyn-s} will be given at the
end of this subsection.

First, we shall present here some algorithm for constructing the sequence
$\Lambda_{1},\Lambda_{2},...,\Lambda_{n}$. By Lemma~\ref{l-on-good-faces}
we can consider only faces $\Lambda$, such that $\overline{\Lambda}=J^{(-)}\times J^{(+)}$.
Algorithm consists of several number of steps and constructs a sequence
$\overline{\Lambda_{1}}$, $\overline{\Lambda_{2}}$, $\ldots$, 
\begin{equation}
\overline{\Lambda}_{p}=J_{p}^{(-)}\times J_{p}^{(+)}=\left\{ (l,k)\,|\, l\in J_{p}^{(-)},\, k\in J_{p}^{(+)}\right\} .\label{eq:L-hat-group}
\end{equation}
 In fact it constructs a sequence $\left\{ (J_{p}^{(-)},J_{p}^{(+)})\right\} _{p=1}^{n}$.
We prefer here to use notation 
\[
(J_{p}^{(-)},J_{p}^{(+)})=T_{p}=\gru{J_{p}^{(-)}}{J_{p}^{(+)}}
\]
 and to call $T_{p}$ a group consisting of particle types listed
in $J_{p}^{(-)},J_{p}^{(+)}$.

Notation $V^{T_{i}}$ has the same meaning as earlier 
\[
V^{T_{i}}\,\,=\,\,\frac{\sum_{l\in J_{i}^{(-)}}v_{l}^{(-)}\rho_{l}^{(-)}+\sum_{k\in J_{i}^{(+)}}v_{k}^{(+)}\rho_{k}^{(+)}}{\sum_{l\in J_{i}^{(-)}}\rho_{l}^{(-)}+\sum_{k\in J_{i}^{(+)}}\rho_{k}^{(+)}}\,.
\]

\emph{Algorithm}: 
\begin{description}
\item [{{{1)}}}] Put $T_{1}=\gru 11$ and find $V^{T_{1}}$ 
\item [{{{2a)}}}] If $V^{T_{1}}<0$, compare $-v_{2}^{(+)}$ and $\left|V^{T_{1}}\right|$.

\begin{itemize}
\item If $-v_{2}^{(+)}>\left|V^{T_{1}}\right|$, then $T_{2}=\gru 1{1,2}$. 
\item If $-v_{2}^{(+)}<\left|V^{T_{1}}\right|$, then $T_{2}=\gru{2,1}1$. 
\end{itemize}
\item [{{{2b)}}}] If $V^{T_{1}}>0$, compare $v_{2}^{(-)}$ and $V^{T_{1}}$.

\begin{itemize}
\item If $v_{2}^{(-)}>V^{T_{1}}$, then $T_{2}=\gru{2,1}1$. 
\item If $v_{2}^{(-)}<V^{T_{1}}$, then $T_{2}=\gru 1{1,2}$. 
\end{itemize}
\end{description}
... Let we have already constructed group 
\[
T_{r-1}=\gru{b,b-1,\ldots,1}{1,\ldots,a-1,a}.
\]
 Find $V^{T_{r-1}}$. If $a<K$ and $b<L$ hold, then apply the following
steps $r$-a) and $r$-b). 
\begin{description}
\item [{{{$r$-a)}}}] If $V^{T_{r-1}}<0$ and $a<K$, compare $-v_{a+1}^{(+)}$
and$\left|V^{T_{r-1}}\right|$.

\begin{itemize}
\item If $-v_{a+1}^{(+)}>\left|V^{T_{r-1}}\right|$, then $T_{r}=\gru{b,\ldots,1}{1,\ldots,a,a+1}$. 
\item If $-v_{a+1}^{(+)}<\left|V^{T_{r-1}}\right|$, then $T_{r}=\gru{b+1,b,\ldots,1}{1,\ldots,a}$. 
\end{itemize}
\item [{{{$r$-b)}}}] If $V^{T_{r-1}}>0$ and $b<L$, we compare $v_{b+1}^{(-)}$
and $V^{T_{r-1}}$.

\begin{itemize}
\item If $v_{b+1}^{(-)}>V^{T_{r-1}}$, then $T_{r}=\gru{b+1,b,\ldots,1}{1,\ldots,a}$. 
\item If $v_{b+1}^{(-)}<V^{T_{r-1}}$, then $T_{r}=\gru{b,\ldots,1}{1,\ldots,a,a+1}$. 
\end{itemize}
\item [{{{$r$-c)}}}] If $a=K$, and $b<L$, we compare $v_{b+1}^{(-)}$
and $V^{T_{r-1}}$.

\begin{itemize}
\item If $v_{b+1}^{(-)}>V^{T_{r-1}}$, then $T_{r}=\gru{b+1,b,\ldots,1}{1,\ldots,K}$. 
\item If $v_{b+1}^{(-)}<V^{T_{r-1}}$, then the algorithm is finished and
the group $T_{r-1}=\gru{b,\ldots,1}{1,\ldots,K}$ is declared to be
the final group $T_{fin}$ of the algorithm. 
\end{itemize}
\item [{{{$r$-d)}}}] If $a<K$, and $b=L$, we compare $v_{a+1}^{(+)}$
and $V^{T_{r-1}}$.

\begin{itemize}
\item If $v_{a+1}^{(+)}<V^{T_{r-1}}$, then $T_{r}=\gru{L,\ldots,1}{1,\ldots,a,a+1}$. 
\item If $v_{a+1}^{(+)}>V^{T_{r-1}}$, then the algorithm is finished and
the group $T_{r-1}=\gru{L,\ldots,1}{1,\ldots,a}$ is declared to be
the final group $T_{fin}$ of the algorithm. 
\end{itemize}
\item [{{{$r$-e)}}}] If $a=K$ and $b=L$, then the algorithm is finished
and the group $T_{r-1}=\gru{L,\ldots,1}{1,\ldots,K}$ is declared
to the final group $T_{fin}$ of the algorithm. 
\end{description}
If the algorithm did not stop at the steps $r$-c), $r$-d) or $r$-e),
then the step $r+1$ should be fulfilled, etc. It is clear that the
algorithm stops after finite number of steps, and as the result we
get a final group $T_{fin}$, which will have one of the following
types
\begin{equation}
\gru{L,\ldots,1}{1,\ldots,K},\qquad\gru{L,\ldots,1}{1,\ldots,K_{1}},\qquad\gru{L_{1},\ldots,1}{1,\ldots,K},\label{eq:fin-gr}
\end{equation}
 where $K_{1}<K$, $L_{1}<L$.

We need not only the final group, corresponding to the face along
which the trajectory escapes to infinity, but also the whole chain
\begin{equation}
T_{1}=\gru 11\rightarrow T_{2}\rightarrow T_{3}\rightarrow\cdots\rightarrow T_{fin}\,.\label{eq:cep-ka}
\end{equation}
 As it follows from the algorithm, this chain is uniquely defined
by the parameters of the model.

Let us remark, that in the algorithm we excluded cases where some
of $V^{T_{r-1}}$ are zero. We will show below (see Remark ~\ref{rem-V-0})
how to modify the algorithm to take into account these cases as well.

The next lemma is needed for the proof of the theorem~\ref{thtrans}.
It is convenient however to give this proof here, as it is essentially
based on the details of the algorithm defined above.

\begin{lemma} \mbox{ } \label{l-vV-1-3} 
\begin{enumerate}
\item If $T_{fin}=\gru{L,\ldots,1}{1,\ldots,K}$, then simultaneously $v_{L}^{(-)}>V^{T_{fin}}$
and~$v_{K}^{(+)}<V^{T_{fin}}$ hold. 
\item If $T_{fin}=\gru{L,\ldots,1}{1,\ldots,K_{1}}$, where $K_{1}<K$,
then $V^{T_{fin}}<0$ and $v_{K}^{(+)}>V^{T_{fin}}$. 
\item If $T_{fin}=\gru{L_{1},\ldots,1}{1,\ldots,K}$, where $L_{1}<L$,
then $V^{T_{fin}}>0$ and $v_{L}^{(-)}<V^{T_{fin}}$. 
\end{enumerate}
\end{lemma}

\emph{Proof of Lemma~\ref{l-vV-1-3}}. In fact, if $T_{fin}=\gru{L,\ldots,1}{1,\ldots,K_{1}}$,
where $K_{1}<K$, then the algorithms stops on some step $r_{0}$-d),
and thus, the condition $v_{K_{1}+1}^{(+)}>V^{T_{fin}}$ will hold.
As $0>v_{K}^{(+)}\geq v_{K_{1}+1}^{(+)}$, then we get the proof of
the part~2 of the lemma. Part~3 is quite similar.

To prove assertion 1 of the lemma consider the face, previous to the
final one. 
\[
T_{fin}=\gru{L,\ldots,1}{1,\ldots,K}.
\]
 Two cases are possible:

\hfill{}$\tfmo=\gru{L,\ldots,1}{1,\ldots,K-1}$~or~$\tfmo=\gru{L-1,\ldots,1}{1,\ldots,K}$.
\hfill{}\null

Consider the case $\tfmo=\gru{L,\ldots,1}{1,\ldots,K-1}$ and the
final fragment of the trajectory in the algorithm:
\[
\gru{L-1,\ldots,1}{1,\ldots,q}\rightarrow\gru{L,\ldots,1}{1,\ldots,q}\rightarrow\cdots\rightarrow\tfmo=\gru{L,\ldots,1}{1,\ldots,K-1}\rightarrow T_{fin}.
\]
 Two cases of the first transition in this chain are possible:

1) $V^{\gru{L-1,\ldots,1}{1,\ldots,q}}<0$ and $v_{q+1}^{(+)}>V^{\gru{L-1,\ldots,1}{1,\ldots,q}}$.

2) $V^{\gru{L-1,\ldots,1}{1,\ldots,q}}>0$ and $v_{L}^{(-)}>V^{\gru{L-1,\ldots,1}{1,\ldots,q}}$.
\\
 In both cases one can claim that 
\begin{equation}
v_{L}^{(-)}>V^{\gru{L,\ldots,1}{1,\ldots,q}}.\label{eq:vL-g-VLq}
\end{equation}
 To prove this consider both cases separately.

Case~1) As $v_{L}^{(-)}>0$, then we have $v_{L}^{(-)}>V^{\gru{L-1,\ldots,1}{1,\ldots,q}}$.
Thus, $v_{L}^{(-)}>V^{\gru{L,\ldots,1}{1,\ldots,q}}$, as $V^{\gru{L,\ldots,1}{1,\ldots,q}}$
is the convex linear combination (CLC%
\footnote{CLC of the numbers $x_{1},\ldots x_{n}$ is $\sum_{i}\alpha_{i}x_{i}$
for some numbers $\alpha_{i}>0$, $i=\overline{1,n}$ such that $\sum_{i}\alpha_{i}=1$.%
}) $v_{L}^{(-)}$ and $V^{\gru{L-1,\ldots,1}{1,\ldots,q}}$.

Case~2) Here we assume $v_{L}^{(-)}>V^{\gru{L-1,\ldots,1}{1,\ldots,q}}$.
From this, as above, we get that $v_{L}^{(-)}>V^{\gru{L,\ldots,1}{1,\ldots,q}}$.

Thus, the inequality (\ref{eq:vL-g-VLq}) is proved. As $V^{\gru{L,\ldots,1}{1,\ldots,K}}$
is CLC of $V^{\gru{L,\ldots,1}{1,\ldots,q}}$ and negative numbers
$v_{q+1}^{(+)}$, $\ldots$, $v_{K}^{(+)}$, then 
\[
V^{\gru{L,\ldots,1}{1,\ldots,K}}<V^{\gru{L,\ldots,1}{1,\ldots,q}}.
\]
 Then we have $V^{\gru{L,\ldots,1}{1,\ldots,K}}<v_{L}^{(-)}$.

\emph{The latter} transition in the chain occurs because $v_{K}^{(+)}<V^{\gru{L,\ldots,1}{1,\ldots,K-1}}$.
Then $v_{K}^{(+)}<V^{\gru{L,\ldots,1}{1,\ldots,K}}$, as $V^{\gru{L,\ldots,1}{1,\ldots,K}}$
is CLC of $V^{\gru{L,\ldots,1}{1,\ldots,K-1}}$ and $v_{K}^{(+)}$.

This gives the proof.

Let $a_{r}$ and $b_{r}$ are such that 
\begin{equation}
T_{r}=\gru{b_{r},\ldots,1}{1,\ldots,a_{r}}.\label{eq:Tr-br-ar}
\end{equation}
 The numbers $a_{r}$ and $b_{r}$ are non-decreasing functions of
$r$. Moreover $a_{r}+b_{r}$ increases by 1 if $r$ increases by~1.
What can be the difference between $T_{r-1}$ and $T_{r}$? There
can be two cases:

Case~$\Pi_{r}$: $\quad$ $a_{r}=a_{r-1}+1$, $b_{r}=b_{r-1}$.

Case~$U_{r}$: $\quad$ $a_{r}=a_{r-1}$, $b_{r}=b_{r-1}+1$.

Remind that the face $B(\Lambda)\in R_{+}^{N}$ is defined by the
set of pairs of indices $\Lambda\subseteq I_{-}\times I_{+}$. Namely,
to each pair $(j,k)\in\Lambda$ corresponds positive coordinates $d_{j,k}>0$
in the definition (\ref{eq:B-L-def1}) of the face $B(\Lambda)$ and
vice-versa. For shortness we say that the face $B(\Lambda)$ consists
of pairs $(j,k)\in\Lambda$.

\noindent \begin{proposition} \label{pr-vect-p-Tr}

\noindent Let the chain~(\ref{eq:cep-ka}) be given and case~$\Pi_{r}$
occurs. For any ergodic face $\Lambda$, not containing the pairs
\begin{equation}
(l,k),\qquad l\in\overline{1,b_{r-1}},\qquad k\in\overline{1,a_{r-1}},\label{eq:uzly}
\end{equation}
 the following holds true: for any pairs as 
\begin{equation}
(b,a_{r}),\quad b\in\overline{1,b_{r-1}}\,,\label{eq:uzly-otr-vvp}
\end{equation}
 belonging to $\Lambda$, the corresponding component of the vector
field is negative : 
\[
v_{b,a_{r}}^{\Lambda}<0\,.
\]
 If the case $U_{r}$ occurs, then for any ergodic face $\Lambda$,
not containing the pairs (\ref{eq:uzly}), the following components
of the vector field are negative 
\[
v_{b_{r},a}^{\Lambda}<0,\quad\quad a\in\overline{1,a_{r-1}}\,,
\]
 under the condition, of course, that $(b_{r},a)\in\Lambda$.

\end{proposition}

\emph{Proof of Proposition \ref{pr-vect-p-Tr}.} Remind the notation
$T_{r}=\gru{b_{r},\ldots,1}{1,\ldots,a_{r}}$. As it was mentioned
above, the connection between $T_{r-1}$ and $T_{r}$ can be of two
kinds~--- $\Pi_{r}$ or $U_{r}$, which we write schematically as
\begin{eqnarray*}
\Pi_{r}: &  & T_{r}=T_{r-1}\cup\gru{\varnothing}{a_{r}}\\
U_{r}: &  & T_{r}=T_{r-1}\cup\gru{b_{r}}{\varnothing}
\end{eqnarray*}
 Consider only the case $\Pi_{r}$, as the case $U_{r}$ is symmetric.
It is necessary to prove that for any ergodic face $\Lambda$, which
does not contain
\[
(l,k),\qquad l\in\overline{1,b_{r-1}},\qquad k\in\overline{1,a_{r-1}},
\]
 for any pairs $(b,a_{r})\in\Lambda$, where $b\in\overline{1,b_{r-1}}$,
the inequality 
\[
v_{b,a_{r}}^{\Lambda}<0\,.
\]
 holds. Thus we mean the faces with 
\[
\overline{\Lambda}=\gru{l_{m},\ldots,l_{r},b_{r-1},\ldots,1}{1,\ldots,a_{r-1},\widehat{a_{r}},k_{r+1},\ldots,k_{n}}.
\]
 For such faces $v_{b,a_{r}}^{\Lambda}=v_{a_{r}}^{(+)}-V^{\overline{\Lambda}}$.

Consider now the case when the set $k_{r+1},\ldots,k_{n}$ is not
empty. As $\overline{\Lambda}$ corresponds to ergodic group of particles,
then by lemma \ref{lln} $v_{k_{r+1}}^{(+)}<V^{\overline{\Lambda}}$.
As $a_{r}<k_{r+1}$, then 
\[
v_{a_{r}}^{(+)}<v_{k_{r+1}}^{(+)}<V^{\overline{\Lambda}}\qquad\Rightarrow\qquad v_{a_{r}}^{(+)}-V^{\overline{\Lambda}}<0.
\]
 The case when the set $k_{r+1},\ldots,k_{n}$ is empty corresponds
to 
\begin{equation}
\overline{\Lambda}=\gru{l_{m},\ldots,l_{r},b_{r-1},\ldots,1}{1,\ldots,a_{r-1}}.\label{eq:dL-lb-1}
\end{equation}
 Case $\Pi_{r}$ includes two possible subcases
\begin{eqnarray}
V^{T_{r-1}}<0, & \quad & v_{a_{r}}^{(+)}<V^{T_{r-1}}\label{eq:A1-vozm}\\
V^{T_{r-1}}>0, & \quad & v_{b_{r-1}+1}^{(-)}<V^{T_{r-1}}\label{eq:A2-nevozm}
\end{eqnarray}
 Consider firstly (\ref{eq:A2-nevozm}). If the set $l_{m},\ldots,l_{r}$
is not empty, then the subcase (\ref{eq:A2-nevozm}) contradicts the
ergodicity assumption for (\ref{eq:dL-lb-1}), thus it is impossible.
If the set $l_{m},\ldots,l_{r}$ is empty, then $\overline{\Lambda}=T_{r-1}$
and the assumption (\ref{eq:A2-nevozm}) means that $V^{\overline{\Lambda}}=V^{T_{r-1}}>0$.
As $v_{a_{r}}^{(+)}<0$, we easily conclude that in this case 
\[
v_{b,a_{r}}^{\Lambda}=v_{a_{r}}^{(+)}-V^{\overline{\Lambda}}<0.
\]
 Consider now (\ref{eq:A1-vozm}). If the set $l_{m},\ldots,l_{r}$
is not empty, then due to the ergodicity of the group (\ref{eq:dL-lb-1}),
we have strict inequality $V^{\overline{\Lambda}}>V^{T_{r-1}}.$ If
the set $l_{m},\ldots,l_{r}$ is empty, then $\overline{\Lambda}=T_{r-1}$
and consequently $V^{\overline{\Lambda}}=V^{T_{r-1}}.$ Finally we
conclude that in the subsituation (\ref{eq:A1-vozm}) always 
\[
V^{\overline{\Lambda}}\geq V^{T_{r-1}}.
\]
 From (\ref{eq:A1-vozm}) we have
\[
v_{a_{r}}^{(+)}<V^{T_{r-1}}\quad\Rightarrow\quad v_{a_{r}}^{(+)}<V^{T_{r-1}}\leq V^{\overline{\Lambda}}
\]
 and it follows that $v_{b,a_{r}}^{\Lambda}=v_{a_{r}}^{(+)}-V^{\overline{\Lambda}}<0$~.

This ends the proof.

\emph{Proof of Proposition~\ref{p-put-dyn-s}}. Assume the above
algorithm produces the chain of groups (\ref{eq:cep-ka}). Let $B(\Lambda_{1})$,
$B(\Lambda_{2})$, $\ldots$, $B(\Lambda_{fin})$ be the faces in
$R_{+}^{N}$, corresponding to the chain $T_{1}$,$T_{2}$, $\ldots$,
$T_{fin}$ via the rule (\ref{eq:L-hat-group}). Denote $F_{1}$,
$F_{2}$, $\ldots$, $F_{fin}$ the closures of these faces in $R_{+}^{N}$.
That is in notation (\ref{eq:Tr-br-ar}) 
\begin{eqnarray*}
F_{i} & = & cl\left(B(\Lambda_{i})\right)=\biggl\{ x\in R_{+}^{N}:\:\,\, x_{j,k}\geq0,\,\,\,(j,k)\notin\left\{ 1,\ldots,b_{r}\right\} \times\left\{ 1,\ldots,a_{r}\right\} ,\\
 &  & \phantom{cl\left(B(\Lambda_{i})\right)=\{\biggl\{ x\in R_{+}^{N}:\:\,\,}\, x_{j,k}=0,\,\,\,(j,k)\in\left\{ 1,\ldots,b_{r}\right\} \times\left\{ 1,\ldots,a_{r}\right\} \,\biggr\}.
\end{eqnarray*}
 It is clear that $F_{1}\supset F_{2}\supset\cdots\supset F_{fin}$,
and moreover $\dim F_{i}>\dim F_{i+1}$. More exactly, $\dim F_{r}-\dim F_{r+1}=b_{r}$
or $a_{r}$ in the case $\Pi_{r}$ or $U_{r}$ correspondingly.

Let $\Gamma_{x}(t)=\left(\dgm_{j,k}(t),\,\,(j,k)\in I_{-}\times I_{+}\right)$
be the coordinate description of the trajectory $\Gamma_{x}$. To
prove that $\Gamma_{x}(t')\in F_{r}$ one should check that $\dgm_{j,k}(t')=0$
for all $(j,k)\in\left\{ 1,\ldots,b_{r}\right\} \times\left\{ 1,\ldots,a_{r}\right\} $.
The trajectory goes along ergodic faces.

1) Maximal ergodic face is $\Lambda_{0}=R_{+}^{N}$. The vector field
$v^{\Lambda_{0}}$ on this face is such that $v_{1,1}^{\Lambda_{0}}=v_{1}^{(+)}-v_{1}^{(-)}<0$.
Note that also for any other ergodic face $\Lambda$, containing the
pair $(1,1)$, the component $v_{1,1}^{\Lambda}$ will also be negative,
as by (\ref{eq:iv1})--(\ref{eq:iv3}) it can take only one of three
following negative values 
\begin{equation}
v_{1}^{(+)}-v_{1}^{(-)},\quad v_{1}^{(+)}-V^{\bar{\Lambda}}\quad\mbox{or}\quad V^{\bar{\Lambda}}-v_{1}^{(-)}.\label{eq:3-neg-val}
\end{equation}
 Thus for any initial point $x$ there is $t_{1}\geq0$ such that
$\dgm_{1,1}(t_{1})=0$, and moreover, $\dgm_{1,1}(t)=0$ $\forall t\geq t_{1}$.

2) Thus $\Gamma_{x}(t_{1})\in F_{1}$. If the case $\Pi_{2}$ occurs,
then we have to show the existence of $t_{2}\geq t_{1}$ such that
$\dgm_{1,2}(t)=0$ $\forall t\geq t_{2}$. If it appeared that $\Gamma_{x}(t_{1})\in F_{2}$,
then just put $t_{2}=t_{1}$. If however $\Gamma_{x}(t_{1})\notin F_{2}$,
that is $\dgm_{1,2}(t_{1})>0$, then $\Gamma_{x}(t_{1})$ belongs
to some ergodic face $\Lambda\ni(1,2)$. By proposition \ref{pr-vect-p-Tr}
$v_{1,2}^{\Lambda}<0$, and thus there is $t_{2}>t_{1}$ such that
$\dgm_{1,2}(t_{2})=0$ (that is $\Gamma_{x}(t_{2})\in F_{2}$). In
future the dynamical system will never quit $F_{2}$. In fact, assume
the contrary. Note that $\Gamma_{x}(t_{2})$ can belong either to
$\Lambda_{2}$, or to its boundary (remind that $\overline{\Lambda_{2}}=\{(1,1),(1,2)\}$
and $cl(\Lambda_{2})=F_{2}$). For the trajectory to quit $F_{2}$
it is necessary that it used some outgoing ergodic face $\Lambda'$.
There are two possibilities to do this. The first possibility is $(1,1)\in\Lambda'$.
But in this case (see~(\ref{eq:3-neg-val})) $v_{1,1}^{\Lambda'}<0$
and we get contradiction with the hypothesis that $\Lambda'$ is an
ergodic outgoing face. The second possibility is $(1,1)\notin\Lambda'$
and $(1,2)\in\Lambda'$. But according to the proposition \ref{pr-vect-p-Tr}
for any such face $v_{1,2}^{\Lambda'}<0$, and thus the dynamical
system cannot quit $F_{2}$ along such face $\Lambda'$, This gives
the contradiction.

If the case $U_{2}$ occurred then, quite similarly, one show existence
of $\, t_{2}\geq t_{1}$ such that $\dgm_{2,1}(t)=0$ $\forall t\geq t_{2}$.

$r$) We can use further the induction, using subsequently proposition
\ref{pr-vect-p-Tr}, to show on the step $r$, that there exists $t_{r}\geq t_{r-1}$
such that for any $t\geq t_{r}$ 
\begin{itemize}
\item $\dgm_{b,a_{r}}(t_{r})=0$ $\forall b\in\overline{1,b_{r-1}}$, if
the case $\Pi_{r}$ holds, 
\item $\dgm_{b_{r},a}(t_{r})=0$ $\forall a\in\overline{1,a_{r-1}}$, if
the case $U_{r}$ holds. 
\end{itemize}
Let us show now that in any case $\Gamma_{x}(t)\in F_{r}$ for all
$t\geq t_{r}$. For concreteness consider only the case $\Pi_{r}$,
that is when 
\begin{eqnarray*}
F_{r-1} & = & \left\{ \, x\in R_{+}^{N}:\,\, x_{i,j}=0\,\:\,\forall(i,j)\in\left\{ b_{r-1},\ldots,1\right\} \times\left\{ 1,\ldots,a_{r-1}\right\} \,\right\} ,\\
F_{r} & = & \left\{ \, x\in R_{+}^{N}:\,\, x_{i,j}=0\,\:\,\forall(i,j)\in\left\{ b_{r-1},\ldots,1\right\} \times\left\{ 1,\ldots,a_{r}\right\} \,\right\} ,\quad\quad a_{r}=a_{r-1}+1.
\end{eqnarray*}
 Assume that the trajectory of the dynamical system $\Gamma_{x}(t)$,
being at time $t=t_{r}$ in $F_{r}$, will leave it at some future
moment. The set $F_{r}$ is a finite union of faces having various
dimensions. One should understand then which outgoing ergodic faces
$\Lambda'$ can be used. Again there are two possibilities..

Case 1: $\Lambda'\bigcap\left\{ b_{r-1},\ldots,1\right\} \times\left\{ 1,\ldots,a_{r-1}\right\} =\varnothing$,
that is $\Lambda'\subset F_{r-1}$. Then there exists $b\in\left\{ b_{r-1},\ldots,1\right\} $
such that $(b,a_{r})\in\Lambda'$ (otherwise $\Lambda'\subset F_{r}$,
which gives the contradiction). By proposition \ref{pr-vect-p-Tr}
we have $v_{b,a_{r}}^{\Lambda'}<0$. This contradicts to the fact
that the face $\Lambda'$ is outgoing.

Case 2: $\Lambda'\bigcap\left\{ b_{r-1},\ldots,1\right\} \times\left\{ 1,\ldots,a_{r-1}\right\} \not=\varnothing$.
Consider 
\[
q=\min\left\{ n:\,\,\Lambda'\bigcap\left\{ b_{n},\ldots,1\right\} \times\left\{ 1,\ldots,a_{n}\right\} \not=\varnothing\,\right\} .
\]
 Assume for definiteness, that on step~$q$ of the algorithm we have
\[
T_{q}=T_{q-1}\cup\gru{b_{q}}{\varnothing}.
\]
 Then there exists such $a\in\left\{ 1,\ldots,a_{q-1}\right\} $,
that $(b_{q},a)\in\Lambda'$. Applying Proposition~\ref{pr-vect-p-Tr},
to $\Lambda'$ we get $v_{b_{q},a}^{\Lambda'}<0$ and come to the
contradiction because $\Lambda'$ is outgoing.

Thus there exists a time moment $t_{fin}>0$ such that for $t\geq t_{fin}$
the trajectory hits the final ergodic face $F_{fin}$, which is the
complement to the final group (\ref{eq:fin-gr}).

Important remark is that the sequence of times
\[
t_{1}\leq t_{2}\leq\cdots\leq\cdots\leq t_{r}\leq\cdots\leq t_{fin}
\]
 depends on the initial point. In particular, for some initial points
some consequent moments $t_{r-1}$ and $t_{r}$ can coincide.

\begin{remark} \emph{\label{rem-V-0} Consider the following modification
of the algorithm: in cases} \textbf{\emph{2a)}} \emph{and $r$-}\textbf{\emph{a)}}
\emph{change the conditions $V^{T_{1}}<0$ and $V^{T_{r-1}}<0$ on
$V^{T_{1}}\leq0$ and $V^{T_{r-1}}\leq0$ correspondingly. All the
rest we leave untouched. It is easy to see that all results of this
section hold after such modification as well. In particular, our study
covers the situation when}%
\footnote{\emph{\,$V^{T_{fin}}$ coincides with the asymptotic boundary velocity
of our system (see Subsection~\ref{sub:proof-th1}).}%
} \emph{$V^{T_{fin}}=0$.}

\end{remark}

From the above it follows that any trajectory $\Gamma_{x}(t)$ reaches
the final face in finite time. To proceed with the proof of Theorem~\ref{l-finite-trans}
we will prove the following lemma.

\begin{lemma} \label{l-finite-trans}For any initial point $x$ the
path $\Gamma_{x}(t)$ has finite number of transitions from one face
to another, until it reaches one of the final faces. In other words
the sequence of faces, passed by the path $\Gamma_{x}(t)$, is finite
and the last element of this sequence is the final face.

\end{lemma}

\emph{Proof of Lemma} \ref{l-finite-trans}. Consider an arbitrary
trajectory $\Gamma_{x}(t)$. Let $\{\Lambda_{i}^{x}\}$ be a sequence
of all faces visited by this trajectory. Denote $\{T_{i}^{x}\}$ the
sequence of the corresponding groups, where $T_{i}^{x}=\overline{\Lambda_{i}^{x}}.$
We want to show that the sequence $\{\Lambda_{i}^{x}\}$ is finite.

Two cases are possible for the transition $\Lambda_{i}^{x}\to\Lambda_{i+1}^{x}$,
or equivalently, for the transition $T_{i}^{x}\to T_{i+1}^{x}$. If
the face $\Lambda_{i}^{x}$ is ergodic, then the group $T_{i+1}^{x}$
is obtained by adding some new particle type to the group $T_{i}^{x}.$
During this transition the dimension of $\Lambda_{i}^{x}$ decreases.
If the face $\Lambda_{i}^{x}$ is non-ergodic, then $\Lambda_{i+1}^{x}$
is the minimal outgoing face, containing $\Lambda_{i}^{x}$ (see lemma
\ref{nonergod}). In the transition $\Lambda_{i}^{x}\to\Lambda_{i+1}^{x}$
from $T_{i}^{x}$ some types are deleted, and the dimension of $\Lambda_{i}^{x}$
increases. Thus, the transition $T_{i}^{x}\to T_{i+1}^{x}$ can occur
with two operations: adding some new type and deleting some types.
The same type can be added and deleted several times. If we could
show that addition and deletion are possible only finite number of
times, that will give finiteness of the sequence $\{\Lambda_{i}^{x}\}$.

Note the following fact. Take for example some $(+)$-type $k$. Then
it can be deleted from the group on some step if and only if on the
previous step we added to the group some $(+)$-type with smaller
number (that is with greater velocity). That is why the type 1, plus
or minus, can be added only once and cannot be deleted. $(+)$-type
2 can be deleted only after adding $(+)$-type 1. Similarly for $(-)$-type
2. That is why type 2, plus or minus, can be added to the group not
more than twice and can be deleted not more than once. One can prove
by induction that any type can be deleted and added not more that
finite number of times.

\emph{Proof of Theorem~\ref{thtrans}}.

Let the chain (\ref{eq:cep-ka}) be the result of the algorithm. Three
cases are possible, defined by simple inequalities between $v_{L}^{(-)}$,
$v_{K}^{(+)}$ and $V^{T_{fin}}$.

$v_{K}^{(+)}<V^{T_{fin}}<v_{L}^{(-)}$. This corresponds to part 1
of lemma \ref{l-vV-1-3}, that is $\Lambda_{fin}=\overline{T_{fin}}=\left\{ 0\right\} .$
Thus (Proposition \ref{p-put-dyn-s}), all trajectories of the dynamical
system $U_{t}$ reach $0$ for finite time and finite number of changes.
Note that from this, using well-known methods (see \cite{VM,FMM}),
one can get alternative proof of ergodicity of $D(t)$, in addition
to the one of theorem \ref{thergod}. The first assertion of theorem
\emph{\ref{thtrans}} is proved .

$V^{T_{fin}}<v_{K}^{(+)}<0$. This case corresponds to part 2 of lemma
\ref{l-vV-1-3}, and thus, 
\[
T_{fin}=\gru{L,\ldots,1}{1,\ldots,K_{1}}
\]
where $K_{1}<K$. From the rules of the algorithm it follows immediately
that $v_{K_{1}+1}^{(+)}>V^{T_{fin}}$, but $v_{K_{1}}^{(+)}<V^{T_{fin}}$.
Thus (see theorem \ref{thergod}), the process $D_{\overline{T_{fin}}}(t)$
is ergodic, and the face $\Lambda_{fin}=\overline{T_{fin}}$ is also
ergodic. Find now the vector $v^{\Lambda_{fin}}$. Note that 
\[
\Lambda_{fin}=\mathcal{L}(L,K_{1})=\{(i,k):\,\, i=1,...,L,\,\, k=K_{1}+1,...,K\}.
\]
 To find components of $v_{i,k}^{\Lambda_{fin}}$ we use the formulas
(\ref{eq:iv1})--(\ref{eq:iv3})
\begin{equation}
v_{i,k}^{\Lambda_{fin}}=v_{k}^{(+)}-V^{\overline{\Lambda_{fin}}}>v_{K_{1}+1}^{(+)}-V^{T_{fin}}>0\qquad\forall\,(i,j)\in\left\{ 1,...,L\right\} \times\left\{ K_{1}+1,...,K\right\} ,\label{eq:v-fin-pol}
\end{equation}
 
\[
v_{i,k}^{\Lambda_{fin}}=0\qquad\forall\,(i,j)\in\left\{ 1,...,L\right\} \times\left\{ 1,...,K_{1}\right\} .
\]
 By Proposition \ref{p-put-dyn-s} any trajectory, in finite time
and after finite number of changes, will reach $\mathcal{L}(L,K_{1})$,
and will move along it with constant speed $v^{\Lambda_{fin}}$, having
strictly positive components (\ref{eq:v-fin-pol}). By standard methods
of \cite{VM,FMM}, we conclude that $D(t)$ is transient. \emph{The
second assertion of theorem \ref{thtrans}} is proved.

$0<v_{L}^{(-)}<V^{T_{fin}}$. This case corresponds to part 3 of lemma
\ref{l-vV-1-3}, and the proof is completely similar to the previous
case. That proves assertion 3 of theorem \emph{\ref{thtrans}}.

The fourth assertion of theorem \emph{\ref{thtrans}} is a corollary
of proposition \ref{p-put-dyn-s} and lemma \ref{l-finite-trans}.

Theorem \emph{\ref{thtrans}} is proved.

\subsection{Proof of theorem 1}

\label{sub:proof-th1}If associated random walk $D(t)$ is ergodic,
then by lemma \ref{lemma1} the speed of the boundary equals $V$
which is defined by (\ref{speed}).

Let the process $D(t)$ be non-ergodic. Then there are two possible
cases: $v_{K}^{(+)}>V$ or $v_{L}^{(-)}<V$. From the previous Subsection
\ref{sub:trans} it follows that any trajectory $\Gamma_{x}(t)$ reaches
the final face in finite time and during this time only finite number
of changing the face occurs.

The following assertion is an obvious analog of the proposition 1.4.3
of \cite{VM}.

\begin{lemma}\label{lln}For any $t\geq0$ and any initial point
$x$
\[
\frac{D_{xM}(tM)}{M}\to\Gamma_{x}(t)
\]
 a.e. as $M\to\infty$.

\end{lemma}

Let $v_{K}^{(+)}>V$. We have proved that any trajectory of the dynamical
system $U^{t}$ reaches the final face $\mathcal{L}(L,K_{1})$, where
the coordinates of the induced vector are positive. By lemma \ref{lln}
the coordinates $d_{q,r}(t)$ of the process $D(t)$, where $q=1,...,L$,
$r=K_{1}+1,...,K,$ grow linearly (a.e.) as $t\in\infty$. In other
words $(+)$-types with numbers $r=K_{1}+1,...,K$ fall behind the
boundary and do not contribute to its velocity. It means that the
boundary velocity is defined only by the particles of types $q=1,...,L$,
$r=1,...,K_{1}$ and are given by formula (\ref{speed}). The case
of $v_{L}^{(-)}<V$ is quite similar.

\section{Appendix }

\subsection{Proof of Lemma \ref{l-on-good-faces}}

Let the face $\Lambda$ be such that $\overline{\Lambda}$ is not
the direct product. Put 
\begin{eqnarray*}
I_{-}^{\overline{\Lambda}} & = & \{i\in I_{-}:\exists k\in I_{+},(i,k)\in\overline{\Lambda}\}\\
I_{+}^{\overline{\Lambda}} & = & \{k\in I_{+}:\exists i\in I_{-},(i,k)\in\overline{\Lambda}\}
\end{eqnarray*}
 Choose an ``appropriate'' face $\Lambda_{0}$ so that $\overline{\Lambda}_{0}=I_{-}^{\overline{\Lambda}}\times I_{+}^{\overline{\Lambda}}$.
To prove the lemma it is sufficient to show 
\[
\mathcal{D}\cap\Lambda=\mathcal{D}\cap\Lambda_{0}
\]
 As $\Lambda\supset\Lambda_{0}$, we always have $\mathcal{D}\cap\Lambda\supset\mathcal{D}\cap\Lambda_{0}$.
Let us prove that $\mathcal{D}\cap\Lambda\subseteq\mathcal{D}\cap\Lambda_{0}$.
Let $(i,k)\in\overline{\Lambda}_{0}$ and $(i,k)\notin\overline{\Lambda}$.
Then there exist $m\in I_{+}$ and $n\in I_{-}$ such that $(i,m)\in\overline{\Lambda}$,
$(n,k)\in\overline{\Lambda}$ and the following equation holds 
\[
d_{i,k}(t)+d_{n,m}(t)=d_{i,m}(t)+d_{n,k}(t)
\]
 Take arbitrary element $d=(d_{j,l})$ of the set $\mathcal{D}\cap\Lambda$.
As its coordinates $d_{i,m}(t)=d_{n,k}(t)=0$, then $d_{i,k}=0$ for
all $(i,k)\in\overline{\Lambda}_{0}$. Thus, $d\in\mathcal{D}\cap\Lambda_{0}$,
and the lemma is proved.

\subsection{Technical lemma \label{sub:Technical-lemma}}

For shortness denote
\[
f(k)=V(I_{-},\{1,...,k\}),\quad g(l)=V(\{1,...,l\},I_{+})
\]

\begin{lemma}\label{tech-lem}We have 
\begin{itemize}
\item $v_{k+1}^{(+)}<f(k+1)\Longleftrightarrow f(k+1)<f(k)$, $k=1,...,K-1$, 
\item $v_{k+1}^{(+)}>f(k+1)\Longleftrightarrow f(k+1)>f(k)$, $k=1,...,K-1$ 
\item $v_{k}^{(+)}>f(k)\Longrightarrow v_{k+1}^{(+)}>f(k+1)$, $k=2,...,K-1$ 
\end{itemize}
Similarly, 
\begin{itemize}
\item $v_{l+1}^{(-)}<g(l+1)\Longleftrightarrow g(l+1)<g(l)$, $l=1,...,L-1$ 
\item $v_{l+1}^{(-)}>g(l+1)\Longleftrightarrow g(l+1)>g(l)$, $l=1,...,L-1$ 
\item $v_{l}^{(-)}<g(l)\Longrightarrow v_{l+1}^{(-)}<g(l+1)$, $l=2,...,L-1$ 
\end{itemize}
\end{lemma}

Proof. We prove the first three items. The others are quite similar.
Using (\ref{speed}) one can check 
\[
f(k+1)=\alpha f(k+1)+\beta f(k+1)=\alpha f(k)+\beta v_{k+1}^{(+)}
\]
 for some $\alpha,\beta>0$ such that $\alpha+\beta=1$. It follows
that
\[
\alpha(f(k+1)-f(k))=\beta(v_{k+1}^{(+)}-f(k+1))
\]
 Thus, $v_{k+1}^{(+)}<f(k+1)\Longleftrightarrow f(k+1)<f(k)$. If
$v_{k}^{(+)}>f(k)$, using $v_{k}^{(+)}<v_{k+1}^{(+)}$, we get 
\[
f(k+1)<\alpha v_{k}^{(+)}+\beta v_{k+1}^{(+)}<\alpha v_{k+1}^{(+)}+\beta v_{k+1}^{(+)}=v_{k+1}^{(+)}
\]

Lemma is proved.

Let $K_{1}$ and $L_{1}$ be defined by (\ref{eq:k1}) and (\ref{eq:l1}).
It follows from the lemma that
\[
f(1)>...>f(K_{1})<f(K_{1}+1)<...<f(K)
\]
 and
\[
g(1)<...<g(L_{1})>g(L_{1}+1)>...>g(L)
\]
 So the minimum of $f(k)$ is reached at point $K_{1}$ and maximum
of $g(l)$ is reached at point $L_{1}$.

\end{document}